\documentclass[draft,11pt]{article}
\usepackage[cp1251]{inputenc}
\usepackage[russian,english]{babel}
\usepackage{amsfonts,amsmath,amsthm,amscd,amssymb,latexsym,
cite,verbatim,texdraw,caption2,pb-diagram}

    \usepackage[usenames]{color}

\theoremstyle{plain}
\newtheorem{remark}{Remark}

\newtheorem{definition}{Definition}

\newtheorem{theorem}{Theorem}
\newtheorem{lemma}[theorem]{Lemma}
\newtheorem{proposition}[theorem]{Proposition}
\newtheorem{corollary}[theorem]{Corollary}
\newtheorem{assumption}[theorem]{Assumption}

\newcommand{\cost}{\text{cost}\medskip}

\newcommand{\z}{\mathbf{z}}

\newcommand{\R}{{\mathbb R} }

\newcommand{\X}{\mathbf{X}}

\newcommand{\K}{\mathsf{K}}
\newcommand{\HK}{\mathcal{H}_{\K}}

\newcommand{\Znu}{\z^{\nu}}
\newcommand{\PZnu}{\mathrm{P}_{\z^{\nu}}}

\newcommand{\F}{\mathcal{F}}

\newcommand{\N}{\mathcal{N}}

\newcommand{\trace}{\text{trace}\medskip}

\numberwithin{equation}{section} \numberwithin{theorem}{section}

\renewcommand{\thetheorem}{\arabic{section}.\arabic{theorem}}
\renewcommand{\theproposition}{\arabic{section}.\arabic{proposition}}
\renewcommand{\thedefinition}{\arabic{section}.\arabic{definition}}
\renewcommand{\thecorollary}{\arabic{section}.\arabic{corollary}}
\renewcommand{\thelemma}{\arabic{section}.\arabic{lemma}}
\renewcommand{\theremark}{\arabic{section}.\arabic{remark}}
\renewcommand{\theexample}{\arabic{section}.\arabic{example}}
\renewcommand{\theequation}{\arabic{section}.\arabic{equation}}

\usepackage{geometry}
\allowdisplaybreaks

\title{On recovering  the Radon-Nikodym derivative under\\ the  big data assumption}

\author{Hanna~L. Myleiko \footnotemark[2]
        \and Sergei~G. Solodky \footnotemark[3]}

\date{}
\begin{document}

\maketitle
\renewcommand{\thefootnote}{\fnsymbol{footnote}}

\footnotetext[2]{ Institute of Mathematics NAS of Ukraine, 3 Tereschenkivska st., Kyiv, Ukraine. Email: hannamyleiko@gmail.com}
\footnotetext[3]{ Institute of Mathematics NAS of Ukraine, 3 Tereschenkivska st., Kyiv, Ukraine; University of Giessen, Department of Mathematics, Giessen, Germany. Email: solodky@imath.kiev.ua}

\begin{abstract}

The present paper is focused on recovering the Radon-Nikodym derivative under the big data assumption.  
To address the above problem, we design an algorithm that is a combination of the Nystr\"om subsampling and the standard  Tikhonov regularization. The convergence rate of the corresponding algorithm is established both in the case when the Radon-Nikodym derivative belongs to RKHS 
and in the case when it does not.  We prove that the proposed approach not only ensures the order of accuracy as algorithms based on the whole sample size, but also allows to achieve subquadratic computational costs in the number of observations.
\end{abstract}

\section{Introduction}
The present study analyzes the implementation of the regularized Nystr\"om subsampling in the context of a numerical approximation of the ratio of two probability density functions, which is usually call the Radon-Nikodym derivative of the corresponding probability measures.
Nowadays, recovering of the Radon-Nikodym derivative is of great interest in  statistical learning since it  can potentially be applied to various tasks such as  transfer learning, covariate shift adaptation, outlier detection, conditional density estimation, etc. Here we may refer to \cite{GizewskiMayer21}, \cite{Shimod00}, \cite{PerBook}, \cite{MylSol03}, \cite{MylSol23}, \cite{SmolSong09},\cite{Hido10}, \cite{Syg10} and the references therein.
At first glance, the simple approach in the density ratio approximation could be performed following the next steps: first, one should estimate  the two density probabilities separately using, for instance, kernel density estimation, and then take the ratio of the obtained estimates.
 However, the algorithmic performance of such an approach is technically more complicated than solving the learning task itself. This disadvantage  is more essential when amount of the involved data is large enough.
In view of the above, a more appropriate approach is one associated with a direct estimation of the density ratio, rather than an estimation of each density separately.

It should be noted that the relevance of the described problem is more significant, the larger amount of data the task deals with. When analyzing such kind of problem the main point is  to reduce storage and  computational costs  arising from big data. The Nystr\"om subsampling is one of the widely used approaches for overcoming these challenges (see, for example, \cite{RudiComRos15}, \cite{KriukPer-Jr16}, \cite{MylSolPer-Jr19}, \cite{LuMathePer-Jr}, \cite{MylSol23}).

In the present study, we are going to employ the regularized Nystr\"om subsampling for recovering the Radon-Nikodym derivative under the big data assumption. 
To the best of our knowledge, up to now, the application of the Nystr\"om family of algorithms to the problem of estimating the Radon-Nikodym derivative was considered only in the context of the domain adaptation with covariate shift (see, e.g., \cite{MylSol03},\cite{MylSol23}). In contrast to the above-mentioned works where the problem of  recovering the Radon-Nikodym derivative was considered under the assumption of high smoothness of the derivative, the present research is devoted to  the case of low smoothness of the derivative. Moreover, when establishing convergence rate of the proposed algorithm, we will take into account both the smoothness of the derivative and the capacity of the space in which it is approximated. In addition, it is worth noting that in the context of the problem under consideration, an interesting aspect for future research  can be an adoptation of the regularized  Nystr\"om subsampling to recover the Radon-Nikodym derivative  with involvement of a more general loss function, such as the one discussed in  \cite{Zel03}, \cite{Grub04}.

The paper is organized as follows. In the next section, we give the strict problem settings and define the Nystr\"om subsampling method. Section 3 contains auxiliary statements and assumptions necessary for further research.  In Section 4 and 5, we obtain error estimates for the regularized Nystr\"om subsampling under the assumptions that the Radon-Nikodym derivative belongs to RKHS, and also does not belong to it, correspondingly. In Section 6, we show that proposed algorithm which is a combination of the Nystr\"om subsampling and the standard  Tikhonov regularization is implemented with subquadratic computational cost on the classes of problems under consideration. Some results allowing a more detailed study of the addressed problem are given in Appendix.

 \section{Problem setting}
 
In the present study, we investigate the problem of  recovering the Radon-Nikodym derivative which can be formulated as follows. 
 Let $p$ and $q$ be two probability measures on a space $\X\subset \R^{d}.$ The information about the measures are only available in the form of samples $\X_{p}=\{x_1,x_2,\ldots,x_N\}$ and $\X_{q}=\{x^{'}_1,x^{'}_2,\ldots,x^{'}_M\}$ drawn independently and identically (i.i.d.) from $p$ and $q$ respectively. Moreover, we assume that there is a function   $\beta\colon \X\rightarrow [0,\infty)$ such that
 $
 dq(x)=\beta(x)dp(x). $ Then $\beta$
  can be viewed as the Radon-Nikodym derivative $\beta=\frac{dq}{dp}.$  
The goal is to approximate  the Radon-Nikodym derivative $\beta$ in two cases 1) $\beta$ belongs to some Reproducing Kernel  Hilbert Space and 2) $\beta$ does not belong to it. Despite  its practical importance, such a case has been studied less in the literature  than the first one. Here we may refer to \cite{LuMathePer-Jr} and the references therein.\\
\subsection{Reproducing Kernel  Hilbert Space}
Let $\HK$ be a Reproducing Kernel  Hilbert Space (RKHS) with symmetric  and positive defined function $\K\colon \X\times \X\rightarrow \mathbb{R}.$ We assume that $\K$ is a continuous and bounded kernel that for any  $x\in \X$ it holds
 \begin{eqnarray}\label{ker}
 \|\K(\cdot,x)\|^{2}_{\HK}=\left<\K(\cdot,x)\textcolor{red}{,} \K(\cdot,x)\right>_{\HK}=\K(x,x)\le\kappa^{2}<\infty.
 \end{eqnarray}
 Let $L_{2,p}$ and  $L_{2,q}$ be  the spaces of the square-integrable function $f\colon \X\rightarrow\mathbb{R} $ with respect to probability measure $p$ and $q$, respectively. 
We also define the canonical embedding operators $J_p\colon \HK\hookrightarrow L_{2,p}$,  $J_q\colon \HK\hookrightarrow L_{2,q}$ and their adjoint operators $J_p^{*}\colon L_{2,p}\rightarrow\HK,$ $J_q^{*}\colon L_{2,q}\rightarrow\HK,$ that  are given  by 
$$
J_p^{*}f(\cdot)=\int_{\X}\K(\cdot,x)f(x)dp(x),\quad\quad\quad
J_q^{*}g(\cdot)=\int_{\X}\K(\cdot,x)g(x)dq(x).
$$

It is known (see, e.g., \cite{NgPer23}, \cite{RudiComRos15}) that in view of (\ref{ker}) for each $\alpha>0$ it holds
$$
\mathcal{N}_{x}(\alpha):=\left\langle \K(\cdot,x),(\alpha I+J_p^{\ast}J_p)^{-1} \K(\cdot,x)\right\rangle_{\HK}=\|(\alpha I+J_p^{\ast}J_p)^{-\frac{1}{2}}\K(\cdot,x)\|^2_{\HK}< \infty.
$$
In the sequel, we define the following quantities
\begin{eqnarray}\label{N_inf}
\mathcal{N}_{\infty}(\alpha):=\sup_{x\in \X} \mathcal{N}_{x}(\alpha)\nonumber
\end{eqnarray}
and
\begin{eqnarray}\label{N_lambda}
\mathcal{N}(\alpha):=\int_{\X}\mathcal{N}_{x}(\alpha) dp(x)=\trace\{(\alpha I+J_p^{\ast}J_p)^{-1} J_p^{\ast}J_p\}.\nonumber
\end{eqnarray}
The function $\mathcal{N}(\alpha)$ measures the capacity of the RKHS $\HK$ in the space $L_{2,p}$ and it is called the effective dimension. 
Further, we formulate the following assumption  which  is  common and not restrictive.  We distinguish two sample operators
$$
S_{\X_{q}}f=(f(x_1^{'}),f(x_2^{'}),\cdots, f(x_M^{'}))\in\mathbb{R}^M,
$$
$$
S_{\X_{p}}f=(f(x_1), f(x_2),\cdots,f(x_N))\in\mathbb{R}^N,
$$
acting from $\mathcal{H}_\K $ to $\mathbb{R}^M$ and $\mathbb{R}^N$, where the norms in  later spaces are 
$M^{-1}$-times and $N^{-1}$-times the standard Euclidian norms, such that the adjoint operators 
 $S^{*}_{\X_{q}}\colon \mathbb{R}^M\rightarrow \mathcal{H}_\K$ and $S^{*}_{\X_p}\colon \mathbb{R}^N\rightarrow \mathcal{H}_\K$ are given as follows
 $$
S^{*}_{\X_{q}}u(\cdot)=\frac{1}{M}\sum_{j=1}^{M}\K(\cdot,x^{'}_{j})u_j,\quad\quad u=(u_1,u_2,\ldots,u_M)\in\mathbb{R}^M,
$$
$$
S^{*}_{\X_{p}}v(\cdot)=\frac{1}{N}\sum_{i=1}^{N}\K(\cdot,x_{i})v_i,\quad\quad v=(v_1,v_2,\ldots,v_N)\in\mathbb{R}^N.
$$

It is easy to see that for any bounded and continuous function $f$ it holds
$$
\int_{\X}\, f(t)\beta(t)dp(t)=\int_{\X}\, f(t)dq(t).
$$
By replacing the function $f(t)$ by $\K(\cdot,t)$, for any $x\in \X$ we get
\begin{equation}\nonumber
\int_{\X}\, \K(x,t)\beta(t)dp(t)=\int_{\X}\, \K(x,t)dq(t),
\end{equation}

\begin{equation}\label{op_eq:1}
J_p^{*}\beta=J_{q}^{*}\mathbf{1},
\end{equation}
where $\mathbf{1}$ is a constant function, that takes the value $1$ everywhere.
Here and in the sequel, we assume that $\mathbf{1}\in\HK$. It should be noted that constant functions may not belong to RKHS $\HK$. This condition is not cumbersome. However, there are $\HK$ such that the constant functions do not belong to them (see, e.g., \cite{Minh10}).

It is known (see, e.g., \cite{LuPer}, \cite{NgPer23}) that  the equation (\ref{op_eq:1}) is ill-posed. Therefore, to  find a solution of (\ref{op_eq:1}) it is necessary to implement algorithms from  the Regularization theory.  

\subsection{Source condition and general regularization scheme}
{\bf{Source condition.}}
Let $\mathcal{H}$ be a Hilbert space  
and  $T\colon \mathcal{H}\rightarrow\mathcal{H},$
be a compact, injective, self-adjoint, and non-negative linear operator. For every $f\in\mathcal{H}$ there is a continuous, strictly increasing function $\phi\colon [0,\|T\|_{\mathcal{H}\rightarrow\mathcal{H}}] \rightarrow \mathbb{R}$, such that  $\phi(0)=0$ and $\phi^2$ is concave. The set of all such functions we denote as $\F.$ 
If $f\in\mathcal{H}$ it can be presented  as 
\begin{eqnarray}\label{source_cond}
f=\phi(T)\mu,\quad\quad \|\mu\|_{\mathcal{H}}\le \varkappa,
\end{eqnarray}
where $\varkappa>0,$ then the expression (\ref{source_cond}) is usually  called  "source condition"$\,$  and $\phi$ is the index function of the source condition (see, e.g., \cite{LuPer}).
Examples of such $\phi$  can be power functions $t^s,\, 0<s\leq\frac{1}{2}$, as well as all less smooth ones. 
In the problem of a numerical representation of the Radon-Nikodym derivative, the source condition was considered in \cite{GizewskiMayer21}, \cite{NgPer23}, \cite{MylSol03}.
\\
{\bf{Regularization scheme.}}
Recall (see \cite{Bak}, \cite{LitPMDiscrStr}) that the most regularization schemes can also be indexed by parameterized function
$g_{\alpha}\colon[0,l]\rightarrow\R$, $l\ge 0,\, \alpha>0$. The only requirements are that there are positive constants $\gamma_{0}, \overline{\gamma}, \tilde{\gamma}$ such that
\begin{equation}\label{eq:g}
\underset{0<t\leq l}{\sup}\,|1-tg_{\alpha}(t)|\leq\gamma_{0},\qquad
\underset{0<t\leq
l}{\sup}\,\sqrt{t}|g_{\alpha}(t)|\leq\frac{\overline{\gamma}}{\sqrt{\alpha}},\qquad
\underset{0<t\leq
l}{\sup}\,|g_{\alpha}(t)|\leq\frac{\tilde{\gamma}}{\alpha} .
\end{equation}

Further important property of the regularization method indexed by $g_{\alpha}$ is its qualification that is the maximum
positive number $p$ for which
\begin{equation}\label{eq:gp}
\underset{0<t\leq l}{\sup}\, t^{p}|1-tg_{\alpha}(t)|\leq
\gamma_{p}\alpha^{p} ,
\end{equation}
where $\gamma_{p}$ does not depend on $\alpha$. 
The following definition  \cite{LuPer, MP2003} shows a relation between the qualification and the source condition.
\begin{definition}
We say that the qualification $p$ covers the index function $\phi$ if the function $t\rightarrow t^p/\phi(t)$ is non-decreasing for 
$t\in (0, l]$.
\end{definition}
The importance of this concept is justified by the following statement.
\begin{proposition}\label{cover}\cite[Proposition 2.7]{LuPer}
Let the regularization method is indexed by $g_{\alpha}(t)$ and has the qualification 
$p$. If this qualification covers the index function 
 $\phi$, then
\begin{equation}\label{cond_qualific} 
\underset{0<t\leq l}{\sup}\,|1-tg_{\alpha}(t)|\phi(t) \leq
\hat{\gamma} \phi(\alpha) ,
\end{equation}
where $\hat{\gamma} = \max\{\gamma_{0}, \gamma_{p}\}$.
\end{proposition}
The proof of this proposition follows directly from (\ref{eq:g}) and  (\ref{eq:gp}).\\

Since the smoothness of such function $f$ is low, then to guarantee  the optimal order of accuracy of its approximation  it is enough to apply a regularization with low qualification $(p=1)$. In our research, as a regularizer we implement the standard Tikhonov method generated by the index function $g_{\alpha}(t)=(t+\alpha)^{-1},\, t,\alpha>0,$ and with the qualification $p=1.$ It should be noted that for the standard Tikhonov method and any index  function $\phi\in\F$ the relations (\ref{cond_qualific}) and
\begin{equation}\label{qualific_root} 
\underset{0<t\leq l}{\sup}\,|1-t g_{\alpha}(t)|\sqrt{t}\phi(t) \leq
\gamma_{*} \sqrt{\alpha}\phi(\alpha)
\end{equation}
hold true.

Hereinafter, we will only consider the function (\ref{source_cond}) with $\phi\in\F.$  

\subsection{Nystr\"om subsampling}
The Nystr\"om type subsampling provides an efficient strategy to conquer the big data challenges. This technique consists of the methods replacing  the entire kernel matrix by a smaller matrix of significantly lower rank, obtained by a random  columns subsampling. 
It is known (see, e.g., \cite{RudiComRos15}) that the Nystr\"om subsampling can be considered as a combination of the standard Tikhonov regularization and a projection scheme on the subset

\begin{eqnarray}\label{HK}
\HK^{\Znu}\colon=\left\{f\colon f=\sum_{i=1}^{|\Znu|}\,c_i\K(\cdot,x_i)+\sum_{j=1}^{|\Znu|}\,c^{'}_j\K(\cdot,x^{'}_j)\right\}, 
\end{eqnarray}
where $|\Znu|\ll \min\{N,M\}.$

Subsequently, for a numerical representation of  the Radon-Nikodym derivative  we  apply the combination of the Nystr\"om  subsampling and the standard Tikhonov regularization. Thus, the approximation to the Radon-Nikodym derivative we will seek as follows 
\begin{eqnarray}\label{eq_beta}
\tilde{\beta}_{M,N,\Znu}^{\alpha_{M,N}}=(\alpha I+\PZnu S^{\ast}_{\X_{p}}S_{\X_{p}}\PZnu)^{-1}\PZnu S^{\ast}_{\X_{q}}S_{\X_{q}}\textbf{1}.
\end{eqnarray}

\section{Auxiliary statements and assumptions}  
In this section we provide some auxiliary statements and assumptions that will be used in the proofs in next sections.
\begin{assumption}\label{source_ker} There is an operator concave index function $\zeta\colon [0,\|J^{*}_{p}J_{p}\|_{\HK\rightarrow\HK}]\rightarrow [0,\infty]$ and $\zeta^{2}$ is covered by qualification $p=1$ such that, for all $x\in \X,$
\begin{equation}\label{ker_source}
\K(\cdot,x)=\zeta(J^{*}_{p}J_{p})\mu_{\K},\qquad \|\mu_{\K}\|_{\HK}\le \overline{\varkappa},
\end{equation}
where $\zeta\in\F$ and $\overline{\varkappa}>0$  does not depend on $x.$
\end{assumption}
Note that the condition (\ref{ker_source}) is the source condition for kernel section $\K(\cdot,x).$ As before, we will consider such $\zeta$ which allows a representation in the power scale $t^r,\, 0<r\leq\frac{1}{2},$ as well as all less smooth ones.  

Here and in the sequel, we adopt the convention that $C$ denotes a generic positive coefficient, which can vary from inequality to inequality and does not depend on the values of $N,M,\alpha,$ and $\delta.$

\begin{lemma}\label{lem:dem}\cite[Lemma 5]{NgPer23}
Under Assumption \ref{source_ker}, it holds
$$
\mathcal{N}_{\infty}(\alpha)\le C\frac{\zeta^{2}(\alpha)}{\alpha}.
$$
\end{lemma}

\begin{lemma}\label{lem:polar}
For operators $J_{p}\colon \HK\hookrightarrow L_{2,p}$ and $J_{p}^{*}\colon  L_{2,p}\rightarrow\HK$ it holds
\begin{eqnarray}\nonumber
\|J_{p}^{*}(\alpha I + J_{p}J_{p}^{*})^{-\frac{1}{2}}\|_{L_{2,p}\rightarrow\HK}
\le 1,
\end{eqnarray}
\begin{eqnarray}\nonumber
\|J_{p}(\alpha I + J_{p}^{*}J_{p})^{-\frac{1}{2}}\|_{\HK\rightarrow L_{2,p}}
\le 1.
\end{eqnarray}
\end{lemma}
The proof of this lemma is given in Appendix \ref{ap:a}.\\\\
Further, we give the following relation for the regularization parameter $\alpha>0,$ sample size $N, $ and the subsample size $|\Znu|.$
For $0<\delta<1$, with probability at least $1-\delta$, we require that
\begin{equation}\label{sample_choice}
|\z^{\nu}|\ge C\mathcal{N}_{\infty}(\alpha)\log\frac{1}{\alpha}\log\frac{1}{\delta}
\end{equation}
and
\begin{equation}\label{par_choice}
\alpha \in \left[C N^{-1}\log\frac{N}{\delta}, \|J^{*}_{p}J_{p}\|_{\HK\rightarrow\HK}\right].
\end{equation}
If $\Znu$ is subsampled according to the plain Nystr\"om approach, then (see, e.g.,\cite[Lemma 6]{RudiComRos15}, \cite[Corollary 1]{MylSolPer-Jr19}) with probability at least $1-\delta$ it holds  
\begin{eqnarray}\label{eq:prob10}
\|(I-\PZnu)(\alpha I+J_p^{*}J_p)^{1/2}\|^{2}_{\HK\rightarrow \HK}\le 3\alpha,
\end{eqnarray}
\begin{eqnarray}\label{eq:prob1}
\|J_p(I-\PZnu)\|^{2}_{\HK\rightarrow L_{2,p}}\le 3\alpha,
\end{eqnarray}
and   for any $\phi\in\F$ (see \cite[Proposition 2]{LitPMDiscrStr})  it holds
\begin{eqnarray}\label{saturat}
\|(I - P_{\z^{\nu}})\phi(J_p^{*}J_p)\|_{\HK\rightarrow\HK}\leq C
\phi(\|(J_p^{*}J_p)^{1/2}(I - P_{\z^{\nu}})\|^{2}_{\HK\rightarrow\HK} ).
\end{eqnarray}
\begin{lemma}\label{lem:add}
For any $\phi\in \F$ it holds
\begin{equation}\label{norm_forHK}
\|(\alpha I+J_p^{*}J_p)^{-1/2}\phi(J_p^{*}J_p)\|_{\HK\rightarrow\HK}\le\frac{1}{\sqrt{\alpha}}\phi(\alpha),
\end{equation}
\begin{equation}\label{norm_forL2}
\|(\alpha I+J_pJ_p^{*})^{-1/2}\phi(J_pJ_p^{*})\|_{L_{2,p}\rightarrow L_{2,p}}\le\frac{1}{\sqrt{\alpha}}\phi(\alpha).
\end{equation}
\end{lemma}
\begin{lemma}\label{lem:add0}
Let $Z\colon \HK\rightarrow\HK$ be a positive self-adjoint operator. For every choice $\Znu$ from the sample $\X_{p}$ we have that
\begin{eqnarray}\label{eq:prob2}
\|(\alpha I+Z)^{\frac{1}{2}}\PZnu(\alpha I+\PZnu Z \PZnu)^{-1}\PZnu(\alpha I+Z)^{\frac{1}{2}}\|_{\HK\rightarrow\HK}\le 1,
\end{eqnarray}
\begin{eqnarray}\label{eq:prob20}
\|(\alpha I+Z)^{\frac{1}{2}}(\alpha I+\PZnu Z\PZnu)^{-1}\PZnu(\alpha I+Z)^{\frac{1}{2}}\|_{\HK\rightarrow\HK}\le 1.
\end{eqnarray}
\end{lemma}
The proof of Lemmas \ref{lem:add} and \ref{lem:add0} are given in Appendix \ref{ap:a}.
\\\\
Following \cite{NgPer23},\cite{LuMathePer}, we introduce the supplemental functions
$$ \mathcal{B}_{N,\alpha}:=\frac{2\kappa}{\sqrt{N}}\left(\frac{\kappa}{\sqrt{N\alpha}}+\sqrt{\mathcal{N}(\alpha)}\right),\quad\quad\quad
  \mathcal{G}(\alpha):=\left(\frac{\mathcal{B}_{N,\alpha}}{\sqrt{\alpha}}\right)^2+1,
 $$
and we will use the auxiliary estimates
\begin{eqnarray}\label{bound:1}
\|(\alpha I + J_{p}^{*}J_{p})^{-\frac{1}{2}}(J_{p}^{*}J_{p}-S_{\X_{p}}^{*}S_{\X_{p}})\|_{\HK\rightarrow\HK}\le\mathcal{B}_{N,\alpha}\log\frac{2}{\delta},
\end{eqnarray}
\begin{eqnarray}\label{bound:2}
\|(\alpha I + J_{p}^{*}J_{p})(\alpha I+S_{\X_{p}}^{*}S_{\X_{p}})^{-1}\|_{\HK\rightarrow\HK}\le 2\left[\left(\frac{\mathcal{B}_{N,\alpha}\log\frac{2}{\delta}}{\sqrt{\alpha}}\right)^{2}+1\right],
\end{eqnarray}
\begin{eqnarray}\label{bound:3}
\|(\alpha I+S_{\X_{p}}^{*}S_{\X_{p}})(\alpha I + J_{p}^{*}J_{p})^{-1}\|_{\HK\rightarrow\HK}\le \frac{\mathcal{B}_{N,\alpha}\log\frac{2}{\delta}}{\sqrt{\alpha}}+1.
\end{eqnarray}

\begin{lemma}\label{lem:add01}
For every choice $\Znu$ from the sample $\X_{p}$ it holds
\begin{eqnarray}\label{eq:prob21}
\|(\alpha I+S_{\X_{p}}^{*}S_{\X_{p}})(\alpha I+\PZnu S_{\X_{p}}^{*}S_{\X_{p}}\PZnu)^{-1}\PZnu\|_{\HK\rightarrow\HK}\le 1+C\left(\frac{\mathcal{B}_{N,\alpha}\log\frac{2}{\delta}}{\sqrt{\alpha}}+1\right)^{\frac{1}{2}}.\nonumber
\end{eqnarray}
\end{lemma}

The proof of Lemma \ref{lem:add01} is given in Appendix \ref{ap:a}.
\begin{proposition}\label{Cordes}\cite[Proposition 4 (Cordes Inequality)]{RudiComRos15}
Let $A,B$ be two positive semidefinite bounded operators on a separable Hilbert space $\mathcal{H}$. Then for all $0\le s\le 1$ it holds
$$
\|A^{s}B^{s}\|_{\mathcal{H}\rightarrow\mathcal{H}}\le\|AB\|_{\mathcal{H}\rightarrow\mathcal{H}}^{s}.
$$
\end{proposition}

\begin{proposition}\label{concentrat}\cite[Proposition 2]{CapVito}
Let $(\Omega, \mathcal{F}, P)$ be a probability space and let $\xi$ be a random variable on $\Omega$ taking value in real separable Hilbert space $H.$ Assume that there are two positive constants $L$ and $\sigma$ such that 
 \begin{eqnarray}\label{concentrat1}
 \mathbb{E}\|\xi-\mathbb{E}\xi\|^{p}_{H}\le \frac{1}{2} p!\sigma^2 L^{p-2}, 
 \end{eqnarray}
for any $p\ge 2.$ Then for all $ l \in \mathbb{N} $ with probability at least $1-\delta$
it holds
 \begin{eqnarray}\label{concentrat2}
 \|\frac{1}{l} \sum_{i=1}^{l}\, \xi(\omega_i)-\mathbb{E}\xi\|_{H}\le 2\left( \frac{L}{l}+\frac{\sigma}{\sqrt{l}}\right)\log\frac{2}{\delta}. \nonumber
 \end{eqnarray}
\end{proposition}
In particular, (\ref{concentrat1}) holds if 
\begin{eqnarray}\label{concentrat3}
\begin{split}
\| \xi(\omega)\|_{H}&\le \frac{L}{2},\quad
\mathbb{E}\|\xi\|_{H}^{2}&\le\sigma^2. \nonumber
\end{split}
\end{eqnarray}

\begin{lemma}\cite[Lemma 1]{NgPer23}\label{lem_eq_prob}
Let $b_0>0$ be such that $|\beta(x)|\le b_0$ for every $x\in \X.$ Then with probability at least $1-\delta$ it holds
 \begin{equation}\label{eq_prob_discr}
 \|(\alpha I + J_{p}^{*}J_{p})^{-\frac{1}{2}}( S_{\X_{p}}^{*}S_{\X_{p}}\beta- S_{\X_{q}}^{*}S_{\X_{q}}\textbf{1})\|_{\HK}\le \left(1+\sqrt{2\log\frac{2}{\delta}}\right)\left(\sqrt{\frac{b_0^2}{N}+\frac{1}{M}}\right)\sqrt{\N_{\infty}(\alpha)}.\nonumber
 \end{equation} 
\end{lemma}
\begin{lemma}\label{lem:b_L2}
For $\beta\in L_2/\HK$, with probability at least $1-\delta$, it holds
 \begin{eqnarray}\label{eq_prob_cont}
\begin{split}
\|(\alpha I + J_{p}^{*}J_{p})^{-\frac{1}{2}}( J_{p}^{*}\beta- S_{\X_{q}}^{*}S_{\X_{q}}\textbf{1})\|_{\HK}
\le C\frac{\sqrt{\mathcal{N}_{\infty}(\alpha)}}{\sqrt{M}}\log^{\frac{1}{2}}\frac{2}{\delta}.\nonumber
\end{split}
\end{eqnarray}
\end{lemma}
The proof of  Lemma \ref{lem:b_L2} is given in Appendix \ref{ap:a}.

\section{Case $\beta\in\HK$}
Recall, that $\beta$ is a solution of the equation (\ref{op_eq:1}). In this Section, we assume that $\beta\in\HK$. For such $\beta$ (see e.g. \cite{KanSuz12}, \cite{SMol20}, \cite{PerBook}),  the equation (\ref{op_eq:1}) can be rewritten as
\begin{equation}\label{eq:HK}
J_{p}^{*}J_p\beta=J_{q}^{*}J_q\mathbf{1}.
\end{equation} 
Note that (\ref{eq:HK}) is ill-posed because the involved operator $J_{p}^{*}J_p$ is compact and its inverse can not be bounded in $\HK.$ 
In this case, it is naturally, to assume that $\beta=\frac{dq}{dp}$ satisfies the source condition 
(\ref{source_cond}) with $T=J_{p}^{*}J_p,$ namely
\begin{eqnarray}\label{source_cond:beta}
\beta=\phi(J_{p}^{*}J_p)\mu_{\beta},
\end{eqnarray}
where $\phi\in\F,$ $\|\mu_{\beta}\|_{\HK}\le \tilde{\varkappa}, \, \tilde{\varkappa}>0.$ 
Recall (see, e.g., \cite{LitPMDiscrStr}) that any $\phi\in \F$ is an operator monotone function, i.e.  for any non-negative self-adjoint operators 
$A,B\colon\HK\rightarrow\HK$ with spectra in $[0,l]$ it holds
\begin{equation} \label{mon_op}
\|\phi(A)-\phi(B)\|_{\HK\rightarrow\HK}\leq C
\phi\left(\|A-B\|_{\HK\rightarrow\HK}\right).
\end{equation}
Now, we are at the point to present main results of this section.
\begin{theorem}\label{H_K}
Assume that in the plain Nystr\"om subsampling the values $|\Znu|$ and $\alpha$ satisfy (\ref{sample_choice}) and (\ref{par_choice}), correspondingly. If $\beta$ satisfies the source condition (\ref{source_cond:beta}), then with probability at least $1-\delta$ it holds 
\begin{eqnarray}\nonumber
\begin{split}
\|\beta-\tilde{\beta}\|_{\HK}&\le C\phi(\alpha)+\left(2\left[\left(\frac{B_{N,\alpha}\log{\frac{2}{\delta}}}{\sqrt{\alpha}}\right)^{2}+1\right]\right)^{1/2}\phi(\alpha)
+C\left(\frac{B_{N,\alpha}\log{\frac{2}{\delta}}}{\sqrt{\alpha}}+1\right)^{1/2}\phi(\alpha)\\
&+\frac{C}{\sqrt{\alpha}}\left(2\left[\left(\frac{B_{N,\alpha}\log{\frac{2}{\delta}}}{\sqrt{\alpha}}\right)^{2}+1\right]\right)^{1/2}\left(\frac{1}{\sqrt{N}}+\frac{1}{\sqrt{M}}\right)\sqrt{\N_{\infty}(\alpha)}\log^{\frac{1}{2}}\frac{2}{\delta},
\end{split}
\end{eqnarray}
\begin{eqnarray}\nonumber
\begin{split}
\|\beta-\tilde{\beta}\|_{L_{2,p}}&\le C\sqrt{\alpha}\phi(\alpha)+\left(2\left[\left(\frac{B_{N,\alpha}\log{\frac{2}{\delta}}}{\sqrt{\alpha}}\right)^{2}+1\right]\right)\sqrt{\alpha}\phi(\alpha)\\
&+C\left(\left(\frac{B_{N,\alpha}\log{\frac{2}{\delta}}}{\sqrt{\alpha}}\right)^2+1\right)^{1/2}\left(\frac{B_{N,\alpha}\log{\frac{2}{\delta}}}{\sqrt{\alpha}}+1\right)^{1/2}\sqrt{\alpha}\phi(\alpha)\\
&+C\left(2\left[\left(\frac{B_{N,\alpha}\log{\frac{2}{\delta}}}{\sqrt{\alpha}}\right)^{2}+1\right]\right)\left(\frac{1}{\sqrt{N}}+\frac{1}{\sqrt{M}}\right)\sqrt{\N_{\infty}(\alpha)},
\end{split}
\end{eqnarray}
where $\tilde{\beta}=\tilde{\beta}_{M,N,\Znu}^{\alpha_{M,N}}$ is  defined by (\ref{eq_beta}).
\end{theorem}
\begin{proof} We split the proof of the theorem into two steps.\\
\text{\bf{Step 1}} \text{\bf{(Estimation in the metric of  $\HK$).}}\\
First, consider the decomposition
\begin{eqnarray}\label{beta_decomp}
\begin{split}
\beta-\tilde{\beta}&=\beta-(\alpha I +\PZnu S_{\X_{p}}^{*}S_{\X_{p}}\PZnu)^{-1}\PZnu S_{\X_{q}}^{*}S_{\X_{q}}\textbf{1}=\omega_0+\omega_1+\omega_2+\omega_3,\nonumber
\end{split}
\end{eqnarray}
where
\begin{eqnarray}\nonumber
\begin{split}
\omega_0&:=( I-\PZnu)\beta;\\
\omega_{1}&:=\PZnu\beta-(\alpha I +\PZnu S_{\X_{p}}^{*}S_{\X_{p}}\PZnu)^{-1}\PZnu S_{\X_{p}}^{*}S_{\X_{p}}\PZnu\beta;\\
\omega_2&:=(\alpha I +\PZnu S_{\X_{p}}^{*}S_{\X_{p}}\PZnu)^{-1}\PZnu S_{\X_{p}}^{*}S_{\X_{p}}\PZnu\beta-(\alpha I +\PZnu S_{\X_{p}}^{*}S_{\X_{p}}\PZnu)^{-1}\PZnu S_{\X_{p}}^{*}S_{\X_{p}}\beta;\\
\omega_3&:=(\alpha I +\PZnu S_{\X_{p}}^{*}S_{\X_{p}}\PZnu)^{-1}\PZnu S_{\X_{p}}^{*}S_{\X_{p}}\beta-(\alpha I +\PZnu S_{\X_{p}}^{*}S_{\X_{p}}\PZnu)^{-1}\PZnu S_{\X_{q}}^{*}S_{\X_{q}}\textbf{1}.
\end{split}
\end{eqnarray}

Now, we estimate the norms of each $\omega_{i},$  $i=\overline{0,3}.$ For $\omega_{0},$ by means of  (\ref{saturat}), we have
\begin{eqnarray}\nonumber
\begin{split}
\|\omega_0\|_{\HK}=\|( I-\PZnu)\beta\|_{\HK}=\|(I-\PZnu)\phi(J_{p}^{*}J_p)\mu_{\beta}\|_{\HK}\le C\phi(\|(J_{p}^{*}J_p)^{\frac{1}{2}}(I-\PZnu)\|^{2}_{\HK\rightarrow\HK}).
\end{split}
\end{eqnarray}
Since the function $\phi(t)$ is covered by the qualification $p=1$, for any $C>1$ we have
\begin{equation}\label{eq:qual}
\frac{t}{\phi(t)} \le \frac{Ct}{\phi(Ct)} \qquad
\Rightarrow \qquad \phi(Ct)\le C \phi(t).
\end{equation}
This together with (\ref{eq:prob1}) implies that
\begin{eqnarray}\label{om:10}
\|\omega_0\|_{\HK}\le C\phi(\alpha).
\end{eqnarray}
Further,
\begin{eqnarray}
\begin{split}
\omega_1&=(I-(\alpha I +\PZnu S_{\X_{p}}^{*}S_{\X_{p}}\PZnu)^{-1}\PZnu S_{\X_{p}}^{*}S_{\X_{p}}\PZnu)\PZnu\beta\\
&=(\alpha I +\PZnu S_{\X_{p}}^{*}S_{\X_{p}}\PZnu)^{-1}\left[\alpha I+ \PZnu S_{\X_{p}}^{*}S_{\X_{p}}\PZnu-\PZnu S_{\X_{p}}^{*}S_{\X_{p}}\PZnu)\right]\PZnu\beta\\
&=\alpha(\alpha I +\PZnu S_{\X_{p}}^{*}S_{\X_{p}}\PZnu)^{-1}\PZnu\beta=
\alpha(\alpha I + S_{\X_{p}}^{*}S_{\X_{p}})^{-\frac{1}{2}}\\
&\times(\alpha I + S_{\X_{p}}^{*}S_{\X_{p}})^{\frac{1}{2}}(\alpha I +\PZnu S_{\X_{p}}^{*}S_{\X_{p}}\PZnu)^{-1}\PZnu
(\alpha I + S_{\X_{p}}^{*}S_{\X_{p}})^{\frac{1}{2}}\\
&\times(\alpha I + S_{\X_{p}}^{*}S_{\X_{p}})^{-\frac{1}{2}}(\alpha I + J_{p}^{*}J_{p})^{\frac{1}{2}}(\alpha I + J_{p}^{*}J_{p})^{-\frac{1}{2}}\beta.
\nonumber
\end{split}
\end{eqnarray}

By (\ref{eq:g}), (\ref{eq:prob20}) with $Z=S_{\X_{p}}^{*}S_{\X_{p}}$, (\ref{bound:2}), (\ref{norm_forHK}) and Proposition \ref{Cordes}, with probability at least $1-\delta$ we obtain
\begin{eqnarray}\label{om:11}
\begin{split}
\|\omega_1\|_{\HK}&\le\alpha\|(\alpha I + S_{\X_{p}}^{*}S_{\X_{p}})^{-\frac{1}{2}}\|_{\HK\rightarrow\HK}\\
&\times\|(\alpha I + S_{\X_{p}}^{*}S_{\X_{p}})^{\frac{1}{2}}(\alpha I +\PZnu S_{\X_{p}}^{*}S_{\X_{p}}\PZnu)^{-1}\PZnu
(\alpha I + S_{\X_{p}}^{*}S_{\X_{p}})^{\frac{1}{2}}\|_{\HK\rightarrow\HK}\\
&\times\|(\alpha I + S_{\X_{p}}^{*}S_{\X_{p}})^{-\frac{1}{2}}(\alpha I + J_{p}^{*}J_{p})^{\frac{1}{2}}\|_{\HK\rightarrow\HK}\|(\alpha I + J_{p}^{*}J_{p})^{-\frac{1}{2}}\phi(J_{p}^{*}J_{p})\mu_{\beta}\|_{\HK}\\
&\le \left(2\left[\left(\frac{B_{N,\alpha}\log{\frac{2}{\delta}}}{\sqrt{\alpha}}\right)^{2}+1\right]\right)^{1/2}\phi(\alpha).
\end{split}
\end{eqnarray}
Next, we are going to bound the norm of $\omega_{2}.$ Recall that
\begin{eqnarray}\nonumber
\omega_2=(\alpha I +\PZnu S_{\X_{p}}^{*}S_{\X_{p}}\PZnu)^{-1}\PZnu S_{\X_{p}}^{*}S_{\X_{p}}(I-\PZnu)\beta.
\end{eqnarray}
Since $\PZnu(I-\PZnu)=0$ we get
\begin{eqnarray}\label{prod_p}
\begin{split}
\PZnu S_{\X_{p}}^{*}S_{\X_{p}}(I-\PZnu)=\PZnu(\alpha I+S_{\X_{p}}^{*}S_{\X_{p}})(I-\PZnu).
\end{split}
\end{eqnarray}
From here, we have
\begin{eqnarray}\nonumber
\begin{split}
\omega_2&=(\alpha I +\PZnu S_{\X_{p}}^{*}S_{\X_{p}}\PZnu)^{-1}\PZnu(\alpha I+S_{\X_{p}}^{*}S_{\X_{p}})(I-\PZnu)\beta\\
&=
(\alpha I + S_{\X_{p}}^{*}S_{\X_{p}})^{-\frac{1}{2}}
(\alpha I + S_{\X_{p}}^{*}S_{\X_{p}})^{\frac{1}{2}}
(\alpha I +\PZnu S_{\X_{p}}^{*}S_{\X_{p}}\PZnu)^{-1}\PZnu
(\alpha I + S_{\X_{p}}^{*}S_{\X_{p}})^{\frac{1}{2}}\\
&\times(\alpha I + S_{\X_{p}}^{*}S_{\X_{p}})^{\frac{1}{2}}(\alpha I + J_{p}^{*}J_{p})^{-\frac{1}{2}}(\alpha I + J_{p}^{*}J_{p})^{\frac{1}{2}}(I-\PZnu)\phi(J_{p}^{*}J_{p})\mu_{\beta},
\end{split}
\end{eqnarray}
then
\begin{eqnarray}\nonumber
\begin{split}
\|\omega_2\|_{\HK}
&\le
\|(\alpha I + S_{\X_{p}}^{*}S_{\X_{p}})^{-\frac{1}{2}}\|_{\HK\rightarrow\HK}\\
&\times\|(\alpha I + S_{\X_{p}}^{*}S_{\X_{p}})^{\frac{1}{2}}
(\alpha I +\PZnu S_{\X_{p}}^{*}S_{\X_{p}}\PZnu)^{-1}\PZnu
(\alpha I + S_{\X_{p}}^{*}S_{\X_{p}})^{\frac{1}{2}}\|_{\HK\rightarrow\HK}\\
&\times\|(\alpha I + S_{\X_{p}}^{*}S_{\X_{p}})^{\frac{1}{2}}(\alpha I + J_{p}^{*}J_{p})^{-\frac{1}{2}}\|_{\HK\rightarrow\HK}\|(\alpha I + J_{p}^{*}J_{p})^{\frac{1}{2}}(I-\PZnu)\|_{\HK\rightarrow\HK}\\
&\times\|(I-\PZnu)\phi(J_{p}^{*}J_{p})\mu_{\beta}\|_{\HK}.
\end{split}
\end{eqnarray}
By means of (\ref{eq:g}), (\ref{eq:prob20}) with $Z=S_{\X_{p}}^{*}S_{\X_{p}}$, (\ref{bound:3}), (\ref{eq:prob10}), Proposition \ref{Cordes}, (\ref{saturat}),  (\ref{eq:prob1}), and (\ref{eq:qual}),
we obtain
\begin{eqnarray}\label{om:12}
\begin{split}
\|\omega_2\|_{\HK}\le\frac{C}{\sqrt{\alpha}}\left(\frac{B_{N,\alpha}\log{\frac{2}{\delta}}}{\sqrt{\alpha}}+1\right)^{1/2}\sqrt{3\alpha}\phi(\alpha)
\le C\left(\frac{B_{N,\alpha}\log{\frac{2}{\delta}}}{\sqrt{\alpha}}+1\right)^{1/2}\phi(\alpha).
\end{split}
\end{eqnarray}

We are at the point to bound the norm of $\omega_{3}.$ We start with the decomposition
\begin{eqnarray}\nonumber
\begin{split}
\omega_3&=(\alpha I +\PZnu S_{\X_{p}}^{*}S_{\X_{p}}\PZnu)^{-1}\PZnu( S_{\X_{p}}^{*}S_{\X_{p}}\beta- S_{\X_{q}}^{*}S_{\X_{q}}\textbf{1})\\
&=
(\alpha I + S_{\X_{p}}^{*}S_{\X_{p}})^{-\frac{1}{2}}
(\alpha I + S_{\X_{p}}^{*}S_{\X_{p}})^{\frac{1}{2}}
(\alpha I +\PZnu S_{\X_{p}}^{*}S_{\X_{p}}\PZnu)^{-1}\PZnu
(\alpha I + S_{\X_{p}}^{*}S_{\X_{p}})^{\frac{1}{2}}\\
&\times(\alpha I + S_{\X_{p}}^{*}S_{\X_{p}})^{-\frac{1}{2}}(\alpha I + J_{p}^{*}J_{p})^{\frac{1}{2}}(\alpha I + J_{p}^{*}J_{p})^{-\frac{1}{2}}( S_{\X_{p}}^{*}S_{\X_{p}}\beta- S_{\X_{q}}^{*}S_{\X_{q}}\textbf{1}),
\end{split}
\end{eqnarray}
then
\begin{eqnarray}\nonumber
\begin{split}
\|\omega_3\|_{\HK}&\le\|(\alpha I + S_{\X_{p}}^{*}S_{\X_{p}})^{-\frac{1}{2}}\|_{\HK\rightarrow\HK}\\
&\times\|(\alpha I + S_{\X_{p}}^{*}S_{\X_{p}})^{\frac{1}{2}}
(\alpha I +\PZnu S_{\X_{p}}^{*}S_{\X_{p}}\PZnu)^{-1}\PZnu
(\alpha I + S_{\X_{p}}^{*}S_{\X_{p}})^{\frac{1}{2}}\|_{\HK\rightarrow\HK}\\
&\times\|(\alpha I + S_{\X_{p}}^{*}S_{\X_{p}})^{-\frac{1}{2}}(\alpha I + J_{p}^{*}J_{p})^{\frac{1}{2}}\|_{\HK\rightarrow\HK}\|(\alpha I + J_{p}^{*}J_{p})^{-\frac{1}{2}}( S_{\X_{p}}^{*}S_{\X_{p}}\beta- S_{\X_{q}}^{*}S_{\X_{q}}\textbf{1})\|_{\HK}.
\end{split}
\end{eqnarray}
By means of (\ref{eq:g}), (\ref{eq:prob20}) with $Z=S_{\X_{p}}^{*}S_{\X_{p}}$, (\ref{bound:2}), Proposition \ref{Cordes}, and Lemma \ref{lem_eq_prob},
with probability at least $1-\delta$ we get 
\begin{eqnarray}\label{om:13}
\|\omega_3\|_{\HK}&\le\frac{C}{\sqrt{\alpha}}\left(2\left[\left(\frac{B_{N,\alpha}\log{\frac{2}{\delta}}}{\sqrt{\alpha}}\right)^{2}+1\right]\right)^{1/2}\left(\frac{1}{\sqrt{N}}+\frac{1}{\sqrt{M}}\right) \sqrt{\N_{\infty}(\alpha)}\log^{\frac{1}{2}}\frac{2}{\delta}.
\end{eqnarray}

Summing up (\ref{om:10}), (\ref{om:11}), (\ref{om:12}), and (\ref{om:13}), we have
\begin{eqnarray}\nonumber
\begin{split}
\|\beta-\tilde{\beta}\|_{\HK}&\le C\phi(\alpha)+\left(2\left[\left(\frac{B_{N,\alpha}\log{\frac{2}{\delta}}}{\sqrt{\alpha}}\right)^{2}+1\right]\right)^{1/2}\phi(\alpha)
+C\left(\frac{B_{N,\alpha}\log{\frac{2}{\delta}}}{\sqrt{\alpha}}+1\right)^{1/2}\phi(\alpha)\\
&+\frac{C}{\sqrt{\alpha}}\left(2\left[\left(\frac{B_{N,\alpha}\log{\frac{2}{\delta}}}{\sqrt{\alpha}}\right)^{2}+1\right]\right)^{1/2}\left(\frac{1}{\sqrt{N}}+\frac{1}{\sqrt{M}}\right)\sqrt{\N_{\infty}(\alpha)}\log^{\frac{1}{2}}\frac{2}{\delta}.
\end{split}
\end{eqnarray}
\\
\text{\bf{Step 2}} \text{\bf{(Estimation in the metric of $L_{2,p}$).}}\\
First, recall (see e.g. \cite[p. 229]{LuPer}) that
$$
\|\beta-\tilde{\beta}\|_{L_{2,p}}=\|(J_p^{*}J_p)^{\frac{1}{2}}(\beta-\tilde{\beta})\|_{\HK}.
$$
Next, similarly to \text{\bf{Step 1}}, we consider the decomposition
\begin{eqnarray}\label{beta_decomp:2}
\begin{split}
(J_p^{*}J_p)^{\frac{1}{2}}(\beta-\tilde{\beta})&=(J_p^{*}J_p)^{\frac{1}{2}}(\beta-(\alpha I +\PZnu S_{\X_{p}}^{*}S_{\X_{p}}\PZnu)^{-1}\PZnu S_{\X_{q}}^{*}S_{\X_{q}}\textbf{1})=\sum_{i=0}^{3}\,\sigma_i\, ,\nonumber
\end{split}
\end{eqnarray}
where
\begin{eqnarray}\nonumber
\begin{split}
\sigma_0&:=(J_p^{*}J_p)^{\frac{1}{2}}( I-\PZnu)\beta;\\
\sigma_{1}&:=(J_p^{*}J_p)^{\frac{1}{2}}\left(\PZnu\beta-(\alpha I +\PZnu S_{\X_{p}}^{*}S_{\X_{p}}\PZnu)^{-1}\PZnu S_{\X_{p}}^{*}S_{\X_{p}}\PZnu\beta\right);\\
\sigma_2&:=(J_p^{*}J_p)^{\frac{1}{2}}\left((\alpha I +\PZnu S_{\X_{p}}^{*}S_{\X_{p}}\PZnu)^{-1}\PZnu S_{\X_{p}}^{*}S_{\X_{p}}\PZnu\beta-(\alpha I +\PZnu S_{\X_{p}}^{*}S_{\X_{p}}\PZnu)^{-1}\PZnu S_{\X_{p}}^{*}S_{\X_{p}}\beta\right);\\
\sigma_3&:=(J_p^{*}J_p)^{\frac{1}{2}}\left((\alpha I +\PZnu S_{\X_{p}}^{*}S_{\X_{p}}\PZnu)^{-1}\PZnu S_{\X_{p}}^{*}S_{\X_{p}}\beta-(\alpha I +\PZnu S_{\X_{p}}^{*}S_{\X_{p}}\PZnu)^{-1}\PZnu S_{\X_{q}}^{*}S_{\X_{q}}\textbf{1}\right).\\
\end{split}
\end{eqnarray}

Now, we estimate the norms of each $\sigma_{i},$  $i=\overline{0,3}.$ 
By means of (\ref{eq:prob1}), (\ref{saturat}), and (\ref{eq:qual}), we obtain
\begin{eqnarray}\label{sg:10}
\begin{split}
\|\sigma_0\|_{\HK}&=\|(J_p^{*}J_p)^{\frac{1}{2}}( I-\PZnu)\beta\|_{\HK}\le\|(J_p^{*}J_p)^{\frac{1}{2}}(I-\PZnu)\|_{\HK\rightarrow\HK}\|(I-\PZnu)\phi(J_{p}^{*}J_p)\mu_{\beta}\|_{\HK}\\
&\le \|(J_p^{*}J_p)^{\frac{1}{2}}(I-\PZnu)\|_{\HK\rightarrow\HK} \phi\left(\|(J_{p}^{*}J_p)^{\frac{1}{2}}(I-\PZnu)\|^{2}_{\HK\rightarrow\HK}\right)\le C\sqrt{\alpha}\phi(\alpha).
\end{split}
\end{eqnarray}
Further,
\begin{eqnarray}
\begin{split}
\sigma_1&=(J_p^{*}J_p)^{\frac{1}{2}}(I-(\alpha I +\PZnu S_{\X_{p}}^{*}S_{\X_{p}}\PZnu)^{-1}\PZnu S_{\X_{p}}^{*}S_{\X_{p}}\PZnu)\PZnu\beta\\
&=(J_p^{*}J_p)^{\frac{1}{2}}(\alpha I +\PZnu S_{\X_{p}}^{*}S_{\X_{p}}\PZnu)^{-1}\left[\alpha I+ \PZnu S_{\X_{p}}^{*}S_{\X_{p}}\PZnu-\PZnu S_{\X_{p}}^{*}S_{\X_{p}}\PZnu)\right]\PZnu\beta\\
&=\alpha (J_p^{*}J_p)^{\frac{1}{2}}(\alpha I +\PZnu S_{\X_{p}}^{*}S_{\X_{p}}\PZnu)^{-1}\PZnu\beta=
\alpha (J_p^{*}J_p)^{\frac{1}{2}}(\alpha I+J_p^{*}J_p)^{-\frac{1}{2}}\\
&\times(\alpha I+J_p^{*}J_p)^{\frac{1}{2}}(\alpha I + S_{\X_{p}}^{*}S_{\X_{p}})^{-\frac{1}{2}}\\
&\times(\alpha I + S_{\X_{p}}^{*}S_{\X_{p}})^{\frac{1}{2}}(\alpha I +\PZnu S_{\X_{p}}^{*}S_{\X_{p}}\PZnu)^{-1}\PZnu
(\alpha I + S_{\X_{p}}^{*}S_{\X_{p}})^{\frac{1}{2}}\\
&\times(\alpha I + S_{\X_{p}}^{*}S_{\X_{p}})^{-\frac{1}{2}}(\alpha I + J_{p}^{*}J_{p})^{\frac{1}{2}}(\alpha I + J_{p}^{*}J_{p})^{-\frac{1}{2}}\beta.
\nonumber
\end{split}
\end{eqnarray}

By Lemma \ref{lem:polar}, (\ref{bound:2}), (\ref{eq:prob20}) with $Z=S_{\X_{p}}^{*}S_{\X_{p}}$, Proposition \ref{Cordes}, and (\ref{norm_forHK}),
 with probability at least $1-\delta$ we have
\begin{eqnarray}\label{sg:11}
\begin{split}
\|\sigma_1\|_{\HK}&\le\alpha\|(J_p^{*}J_p)^{\frac{1}{2}}(\alpha I+J_p^{*}J_p)^{-\frac{1}{2}}\|_{\HK\rightarrow\HK}\\
&\times\|(\alpha I+J_p^{*}J_p)^{\frac{1}{2}}(\alpha I + S_{\X_{p}}^{*}S_{\X_{p}})^{-\frac{1}{2}}\|_{\HK\rightarrow\HK}\\
&\times\|(\alpha I + S_{\X_{p}}^{*}S_{\X_{p}})^{\frac{1}{2}}(\alpha I +\PZnu S_{\X_{p}}^{*}S_{\X_{p}}\PZnu)^{-1}\PZnu
(\alpha I + S_{\X_{p}}^{*}S_{\X_{p}})^{\frac{1}{2}}\|_{\HK\rightarrow\HK}\\
&\times\|(\alpha I + S_{\X_{p}}^{*}S_{\X_{p}})^{-\frac{1}{2}}(\alpha I + J_{p}^{*}J_{p})^{\frac{1}{2}}\|_{\HK\rightarrow\HK}\|(\alpha I + J_{p}^{*}J_{p})^{-\frac{1}{2}}\phi(J_{p}^{*}J_{p})\mu_{\beta}\|_{\HK}\\
&\le \sqrt{\alpha}\left(2\left[\left(\frac{B_{N,\alpha}\log{\frac{2}{\delta}}}{\sqrt{\alpha}}\right)^{2}+1\right]\right)\phi(\alpha).
\end{split}
\end{eqnarray}
Next, we are going to bound the norm of $\sigma_{2}.$ Recall that
\begin{eqnarray}\nonumber
\sigma_2=(J_p^{*}J_p)^{\frac{1}{2}}(\alpha I +\PZnu S_{\X_{p}}^{*}S_{\X_{p}}\PZnu)^{-1}\PZnu S_{\X_{p}}^{*}S_{\X_{p}}(I-\PZnu)\beta.
\end{eqnarray}
Using (\ref{prod_p}), we have
\begin{eqnarray}\nonumber
\begin{split}
\sigma_2&=(J_p^{*}J_p)^{\frac{1}{2}}(\alpha I +\PZnu S_{\X_{p}}^{*}S_{\X_{p}}\PZnu)^{-1}\PZnu(\alpha I+S_{\X_{p}}^{*}S_{\X_{p}})(I-\PZnu)\beta\\
&=(J_p^{*}J_p)^{\frac{1}{2}}(\alpha I+J_p^{*}J_p)^{-\frac{1}{2}}
(\alpha I+J_p^{*}J_p)^{\frac{1}{2}}
(\alpha I + S_{\X_{p}}^{*}S_{\X_{p}})^{-\frac{1}{2}}\\
&\times(\alpha I + S_{\X_{p}}^{*}S_{\X_{p}})^{\frac{1}{2}}
(\alpha I +\PZnu S_{\X_{p}}^{*}S_{\X_{p}}\PZnu)^{-1}\PZnu
(\alpha I + S_{\X_{p}}^{*}S_{\X_{p}})^{\frac{1}{2}}\\
&\times(\alpha I + S_{\X_{p}}^{*}S_{\X_{p}})^{\frac{1}{2}}(\alpha I + J_{p}^{*}J_{p})^{-\frac{1}{2}}(\alpha I + J_{p}^{*}J_{p})^{\frac{1}{2}}(I-\PZnu)\phi(J_{p}^{*}J_{p})\mu_{\beta},
\end{split}
\end{eqnarray}
then
\begin{eqnarray}\nonumber
\begin{split}
\|\sigma_2\|_{\HK}
&\le
\|(J_p^{*}J_p)^{\frac{1}{2}}(\alpha I+J_p^{*}J_p)^{-\frac{1}{2}}\|_{\HK\rightarrow\HK}
\|(\alpha I+J_p^{*}J_p)^{\frac{1}{2}}(\alpha I + S_{\X_{p}}^{*}S_{\X_{p}})^{-\frac{1}{2}}\|_{\HK\rightarrow\HK}\\
&\times\|(\alpha I + S_{\X_{p}}^{*}S_{\X_{p}})^{\frac{1}{2}}
(\alpha I +\PZnu S_{\X_{p}}^{*}S_{\X_{p}}\PZnu)^{-1}\PZnu
(\alpha I + S_{\X_{p}}^{*}S_{\X_{p}})^{\frac{1}{2}}\|_{\HK\rightarrow\HK}\\
&\times\|(\alpha I + S_{\X_{p}}^{*}S_{\X_{p}})^{\frac{1}{2}}(\alpha I + J_{p}^{*}J_{p})^{-\frac{1}{2}}\|_{\HK\rightarrow\HK}\|(\alpha I + J_{p}^{*}J_{p})^{\frac{1}{2}}(I-\PZnu)\|_{\HK\rightarrow\HK}\\
&\times\|(I-\PZnu)\phi(J_{p}^{*}J_{p})\mu_{\beta}\|_{\HK}.
\end{split}
\end{eqnarray}
Applying  Lemma \ref{lem:polar}, (\ref{bound:2}), (\ref{eq:prob20}) with $Z=S_{\X_{p}}^{*}S_{\X_{p}}$, (\ref{bound:3}), (\ref{eq:prob10}), (\ref{saturat}),  Proposition \ref{Cordes}, (\ref{eq:prob1}), and (\ref{eq:qual}), with probability at least $1-\delta$
we obtain
\begin{eqnarray}\label{sg:12}
\begin{split}
\|\sigma_2\|_{\HK}\le C\left(\left(\frac{B_{N,\alpha}\log{\frac{2}{\delta}}}{\sqrt{\alpha}}\right)^2+1\right)^{1/2}\left(\frac{B_{N,\alpha}\log{\frac{2}{\delta}}}{\sqrt{\alpha}}+1\right)^{1/2}\sqrt{\alpha}\phi(\alpha).
\end{split}
\end{eqnarray}

We are at the point to bound the norm of $\sigma_{3}.$ We start with the decomposition
\begin{eqnarray}\nonumber
\begin{split}
\sigma_3&=(J_p^{*}J_p)^{\frac{1}{2}}(\alpha I +\PZnu S_{\X_{p}}^{*}S_{\X_{p}}\PZnu)^{-1}\PZnu( S_{\X_{p}}^{*}S_{\X_{p}}\beta- S_{\X_{q}}^{*}S_{\X_{q}}\textbf{1})\\
&=(J_p^{*}J_p)^{\frac{1}{2}}(\alpha I+J_p^{*}J_p)^{-\frac{1}{2}}
(\alpha I+J_p^{*}J_p)^{\frac{1}{2}}
(\alpha I + S_{\X_{p}}^{*}S_{\X_{p}})^{-\frac{1}{2}}\\
&\times(\alpha I + S_{\X_{p}}^{*}S_{\X_{p}})^{\frac{1}{2}}
(\alpha I +\PZnu S_{\X_{p}}^{*}S_{\X_{p}}\PZnu)^{-1}\PZnu
(\alpha I + S_{\X_{p}}^{*}S_{\X_{p}})^{\frac{1}{2}}\\
&\times(\alpha I + S_{\X_{p}}^{*}S_{\X_{p}})^{-\frac{1}{2}}(\alpha I + J_{p}^{*}J_{p})^{\frac{1}{2}}(\alpha I + J_{p}^{*}J_{p})^{-\frac{1}{2}}( S_{\X_{p}}^{*}S_{\X_{p}}\beta- S_{\X_{q}}^{*}S_{\X_{q}}\textbf{1}),
\end{split}
\end{eqnarray}
then
\begin{eqnarray}\nonumber
\begin{split}
\|\sigma_3\|_{\HK}&\le\|(J_p^{*}J_p)^{\frac{1}{2}}(\alpha I+J_p^{*}J_p)^{-\frac{1}{2}}\|_{\HK\rightarrow \HK}
\|(\alpha I+J_p^{*}J_p)^{\frac{1}{2}}(\alpha I + S_{\X_{p}}^{*}S_{\X_{p}})^{-\frac{1}{2}}\|_{\HK\rightarrow\HK}\\
&\times\|(\alpha I + S_{\X_{p}}^{*}S_{\X_{p}})^{\frac{1}{2}}
(\alpha I +\PZnu S_{\X_{p}}^{*}S_{\X_{p}}\PZnu)^{-1}\PZnu
(\alpha I + S_{\X_{p}}^{*}S_{\X_{p}})^{\frac{1}{2}}\|_{\HK\rightarrow\HK}\\
&\times\|(\alpha I + S_{\X_{p}}^{*}S_{\X_{p}})^{-\frac{1}{2}}(\alpha I + J_{p}^{*}J_{p})^{\frac{1}{2}}\|_{\HK\rightarrow\HK}\|(\alpha I + J_{p}^{*}J_{p})^{-\frac{1}{2}}( S_{\X_{p}}^{*}S_{\X_{p}}\beta- S_{\X_{q}}^{*}S_{\X_{q}}\textbf{1})\|_{\HK}.
\end{split}
\end{eqnarray}
By means of  Lemma \ref{lem:polar}, (\ref{eq:prob20}) with $Z=S_{\X_{p}}^{*}S_{\X_{p}}$, (\ref{bound:2}), Proposition \ref{Cordes}, and  Lemma \ref{lem_eq_prob},
with probability at least $1-\delta$ we get 
\begin{eqnarray}\label{sg:13}
\begin{split}
\|\sigma_3\|_{\HK}&\le C\left(2\left[\left(\frac{B_{N,\alpha}\log{\frac{2}{\delta}}}{\sqrt{\alpha}}\right)^{2}+1\right]\right)\left(\frac{1}{\sqrt{N}}+\frac{1}{\sqrt{M}}\right)\sqrt{\N_{\infty}(\alpha)}.
\end{split}
\end{eqnarray}

Summing up (\ref{sg:10}), (\ref{sg:11}), (\ref{sg:12}), and (\ref{sg:13}), we have
\begin{eqnarray}\nonumber
\begin{split}
\|\beta-\tilde{\beta}\|_{L_{2,p}}&\le C\sqrt{\alpha}\phi(\alpha)+\left(2\left[\left(\frac{B_{N,\alpha}\log{\frac{2}{\delta}}}{\sqrt{\alpha}}\right)^{2}+1\right]\right)\sqrt{\alpha}\phi(\alpha)\\
&+C\left(\left(\frac{B_{N,\alpha}\log{\frac{2}{\delta}}}{\sqrt{\alpha}}\right)^2+1\right)^{1/2}\left(\frac{B_{N,\alpha}\log{\frac{2}{\delta}}}{\sqrt{\alpha}}+1\right)^{1/2}\sqrt{\alpha}\phi(\alpha)\\
&+C\left(2\left[\left(\frac{B_{N,\alpha}\log{\frac{2}{\delta}}}{\sqrt{\alpha}}\right)^{2}+1\right]\right)\left(\frac{1}{\sqrt{N}}+\frac{1}{\sqrt{M}}\right)\sqrt{\N_{\infty}(\alpha)}.
\end{split}
\end{eqnarray}
\end{proof}

In the sequel, we will use the following statement
\begin{lemma}\cite[Lemma 4.6]{LuMathePer}\label{lem:LuPer}
There exists $\alpha_{*}$ such that $\frac{\mathcal{N}(\alpha_{*})}{\alpha_{*}}=N.$ 
For 
$\alpha_{*}\le\alpha\le\kappa$
 there holds
\begin{equation}\label{est:B}
\mathcal{B}_{N,\alpha}\le\frac{2\kappa}{\sqrt{N}}(\sqrt{2}\kappa+\sqrt{\mathcal{N}(\alpha)}).
\end{equation}
This yields
\begin{equation}\label{est:G}
\mathcal{G}(\alpha)\le1+(4\kappa^2+2\kappa)^{2}\nonumber
\end{equation}
and also 
\begin{equation}\label{est:B2}
\mathcal{B}_{N,\alpha}(\mathcal{B}_{N,\alpha}+\sqrt{\alpha})\le(1+4\kappa)^{4}\min\left\{\alpha,\sqrt{\frac{\kappa}{N}}\right\}.\nonumber
\end{equation}
\end{lemma}

\begin{theorem}\label{ker_HK}
Let $\K$ satisfies Assumption \ref{source_ker} and  $\alpha\ge\alpha_{*}$. Then under the assumption of Theorem \ref{H_K} and if $\alpha_{*}$ obeys (\ref{par_choice})  with probability at least $1-\delta$ it holds
\begin{eqnarray}\label{err_HK}
\begin{split}
\|\beta-\tilde{\beta}\|_{\HK}&\le C\Big(\phi(\alpha)+\left(\frac{1}{\sqrt{N}}+\frac{1}{\sqrt{M}}\right)\frac{\zeta(\alpha)}{\alpha}\Big)\log^{2}\frac{2}{\delta},
\end{split}
\end{eqnarray}
\begin{eqnarray}\label{err_L2}
\begin{split}
\|\beta-\tilde{\beta}\|_{L_{2,p}}&\le C\sqrt{\alpha}\Big(\phi(\alpha)+\left(\frac{1}{\sqrt{N}}+\frac{1}{\sqrt{M}}\right)\frac{\zeta(\alpha)}{\alpha}\Big)\log^{2}\frac{2}{\delta},\nonumber
\end{split}
\end{eqnarray}
where $\zeta$ is defined by (\ref{ker_source}).
\end{theorem}
\begin{proof}
From Theorem \ref{H_K} and  Lemma \ref{lem:LuPer} it follows
\begin{eqnarray}\nonumber
\begin{split}
\|\beta-\tilde{\beta}\|_{\HK}&\le C\phi(\alpha)+\frac{C}{\sqrt{\alpha}}\log^{2}\frac{2}{\delta}\left(\frac{1}{\sqrt{N}}+\frac{1}{\sqrt{M}}\right)\sqrt{\N_{\infty}(\alpha)}
\end{split}
\end{eqnarray}
and
\begin{eqnarray}\nonumber
\begin{split}
\|\beta-\tilde{\beta}\|_{L_{2,p}}&\le C\sqrt{\alpha}\phi(\alpha)+C\log^{2}\frac{2}{\delta}\left(\frac{1}{\sqrt{N}}+\frac{1}{\sqrt{M}}\right)\sqrt{\N_{\infty}(\alpha)}.
\end{split}
\end{eqnarray}
Further, applying Lemma \ref{lem:dem} we get the statement of Theorem. 
\end{proof}

\begin{corollary}\label{col_HK}
Denote $\theta_{\phi,\zeta}(t)=\frac{t\phi(t)}{\zeta(t)}$ and $\alpha=\alpha_{N,M}:=\theta_{\phi,\zeta}^{-1}\left(\frac{1}{\sqrt{N}}+\frac{1}{\sqrt{M}}\right),$ then with probability at least $1-\delta$ it holds
\begin{equation}\label{err_HK1}
\|\beta-\tilde{\beta}\|_{\HK}\le C\phi\left(\theta_{\phi,\zeta}^{-1}\left(\frac{1}{\sqrt{N}}+\frac{1}{\sqrt{M}}\right)\right)\log^2\frac{2}{\delta},
\end{equation}
\begin{equation}\label{err_L2_1}
\|\beta-\tilde{\beta}\|_{L_{2,p}}\le C\sqrt{\theta_{\phi,\zeta}^{-1}\left(\frac{1}{\sqrt{N}}+\frac{1}{\sqrt{M}}\right)}\phi\left(\theta_{\phi,\zeta}^{-1}\left(\frac{1}{\sqrt{N}}+\frac{1}{\sqrt{M}}\right)\right)\log^2\frac{2}{\delta}.
\end{equation}
Furthermore, if $\beta$ meets the source condition (\ref{source_cond}) with $\phi(t)=t^{s},\, s\in(0,\frac{1}{2} ],$ and $\K$ satisfies Assumption \ref{source_ker} with $\zeta(t)=t^{r},\, r\in(0,\frac{1}{2}],$ then 
\begin{eqnarray}\label{t_sr:hk1}
\|\beta-\tilde{\beta}\|_{\HK}\le C\left(\frac{1}{\sqrt{N}}+\frac{1}{\sqrt{M}}\right)^{\frac{s}{s+1-r}}\log^{2}\frac{2}{\delta},
\end{eqnarray}
\begin{equation}\label{t_sr:l21}
\|\beta-\tilde{\beta}\|_{L_{2,p}}\le C\left(\frac{1}{\sqrt{N}}+\frac{1}{\sqrt{M}}\right)^{\frac{1+2s}{2(s+1-r)}}\log^{2}\frac{2}{\delta}.
\end{equation}
\end{corollary}
\begin{proof}
Let
$$
\phi(\alpha)=\left(\frac{1}{\sqrt{N}}+\frac{1}{\sqrt{M}}\right)\frac{\zeta(\alpha)}{\alpha}.
$$
Then
$$
\frac{1}{\sqrt{N}}+\frac{1}{\sqrt{M}}=\frac{\alpha\phi(\alpha)}{\zeta(\alpha)} \quad \Rightarrow \quad \alpha=\theta_{\phi,\zeta}^{-1}\left(\frac{1}{\sqrt{N}}+\frac{1}{\sqrt{M}}\right),$$
with $\theta_{\phi,\zeta}(t)=\frac{t\phi(t)}{\zeta(t)}$.

Thus, from Theorem \ref{ker_HK} it follows
$$
\|\beta-\tilde{\beta}\|_{\HK}\le C\phi\left(\theta_{\phi,\zeta}^{-1}\left(\frac{1}{\sqrt{N}}+\frac{1}{\sqrt{M}}\right)\right)\log^2\frac{2}{\delta}.
$$
Using the same arguments as for the establishing (\ref{err_HK1}), one can obtain (\ref{err_L2_1}).  
Further, substituting $\phi(t)=t^{s}$ and $\zeta(t)=t^{r}$ into (\ref{err_HK1}) and (\ref{err_L2_1}), we derive to  (\ref{t_sr:hk1}) and (\ref{t_sr:l21}), correspondingly.
\end{proof}
Note that the both bounds in Theorem \ref{ker_HK} are valid for all $\alpha>\alpha_{*}.$
Applying Lemmas \ref{lem:LuPer} and \ref{lem:dem} it is easy to check that $\alpha=\alpha_{N,M}$ also satisfies the above inequality. 

\begin{remark}
Earlier, the problem of approximating the Radon-Nikodym derivative in the metric of $\HK$ was studied in the works
\cite{GizewskiMayer21},\cite{NgPer23},\cite{MylSol03}. Unlike  \cite{GizewskiMayer21}, \cite{NgPer23}, in \cite{MylSol03} this problem was considered in the  Big Data setting as well as in our present research. In addition, in \cite{MylSol03} it was assumed that the exact solution satisfies the source condition (\ref{source_cond:beta}) with smoother than we assume index functions, namely, functions representable in a power scale in the form $t^s , 1\le s\le \frac{3}{2}.$
Note, the accuracy estimate in (\ref{err_HK}) coincides, up to a constant, with the estimates from \cite{GizewskiMayer21},\cite{NgPer23}, \cite{MylSol03}. As for the estimate in the metric of $L_{2,p}$, the similar research was conducted in \cite{KanSuz12}. In the present analysis, we use the source condition, which allows us to obtain more accurate error estimate than in \cite{KanSuz12}. Moreover, we continue and extend the previous studies to the  new classes of problems. 
\end{remark}

\section{Case $\beta\in L_{2,p}/\HK$}
In this section we consider the case when $\beta\in L_{2,p}/\HK$, i.e., $\beta$ does not belong to RKHS $\HK$. 
Here we assume that  $\beta=\frac{dq}{dp}$ satisfies the source condition 
(\ref{source_cond}) with $T=J_{p}J^{*}_p,$ namely
\begin{eqnarray}\label{source_cond:beta2}
\beta=\phi(J_{p}J^{*}_p)\mu_{\beta},
\end{eqnarray}
where $\phi\in \F, \quad \|\mu_{\beta}\|_{L_{2,p}}\le\hat{\varkappa},\,\hat{\varkappa}>0.$ From the definition of $\F$ it follows that $\frac{\sqrt{t}}{\phi(t)}$ is non-decreasing, which means that $\phi(t)$ increases not faster than $\sqrt{t}$.

\begin{theorem}\label{L_2}
Assume that in the plain Nystr\"om subsampling the values $|\Znu|$ and $\alpha$ satisfy (\ref{sample_choice}), (\ref{par_choice}), correspondingly. If $\beta$ satisfies the source condition (\ref{source_cond:beta2}), then with probability at least $1-\delta$ it holds 
 \begin{eqnarray}\nonumber
 \begin{split}
\|J_{p}^{*}(\beta&-J_{p}\tilde{\beta})\|_{\HK}\le C\sqrt{\alpha}\phi(\alpha)\\
&+ 
C\left(1+C\left(\frac{\mathcal{B}_{N,\alpha}\log\frac{2}{\delta}}{\sqrt{\alpha}}+1\right)^{\frac{1}{2}}\right)\left[\left(\frac{B_{N,\alpha}\log\frac{2}{\delta}}{\sqrt{\alpha}}\right)^{2}+1\right]\left(\frac{1}{\sqrt{M}}+
\frac{1}{\sqrt{M}}\right)\sqrt{\mathcal{N}_{\infty}(\alpha)}\log^{\frac{1}{2}}\frac{2}{\delta},
\end{split}
\end{eqnarray}
\begin{eqnarray}\nonumber
\begin{split}
\|\beta-J_{p}\tilde{\beta}\|_{L_{2,p}}&\le C\phi(\alpha)
+ C\left[\left(\frac{B_{N,\alpha}\log{\frac{2}{\delta}}}{\sqrt{\alpha}}\right)^2+1\right]\frac{B_{N,\alpha}}{\sqrt{\alpha}}\phi(\alpha)\log^{\frac{1}{2}}{\frac{2}{\delta}}\\
&+ C\left[\left(\frac{B_{N,\alpha}\log\frac{2}{\delta}}{\sqrt{\alpha}}\right)^{2}+1\right]
\frac{\sqrt{\mathcal{N}_{\infty}(\alpha)}}{\sqrt{M}}\log^{\frac{1}{2}}\frac{2}{\delta},
\end{split}
\end{eqnarray}
where $\tilde{\beta}=\tilde{\beta}_{M,N,\Znu}^{\alpha_{M,N}}$ is defined by (\ref{eq_beta}).
\end{theorem}
The proof of  Theorem \ref{L_2} is deferred to Appendix \ref{ap_b}.

\begin{theorem}\label{ker_L2}
Let $\K$ satisfies Assumption \ref{source_ker} and $\alpha\ge\alpha_{*}$. Then under the assumption of Theorem \ref{L_2} and if $\alpha_{*}$ obeys (\ref{par_choice}) with probability at least $1-\delta$ it holds\\
\begin{eqnarray}\label{n_hk}
\begin{split}
\|J_{p}^{*}(\beta-J_{p}\tilde{\beta})\|_{\HK}\le C\Big(\sqrt{\alpha}\phi(\alpha)+
\left(\frac{1}{\sqrt{N}}+
\frac{1}{\sqrt{M}}\right)\frac{\zeta(\alpha)}{\sqrt{\alpha}}\Big)\log^{2}\frac{2}{\delta},\nonumber
\end{split}
\end{eqnarray}
\begin{eqnarray}\label{n_l2}
\begin{split}
\|\beta-J_{p}\tilde{\beta}\|_{L_{2,p}}\le C\Big(\phi(\alpha)+
\frac{1}{\sqrt{M}}\frac{\zeta(\alpha)}{\sqrt{\alpha}}\Big)\log^{2}\frac{2}{\delta},\nonumber
\end{split}
\end{eqnarray}
where $\zeta$ is defined by (\ref{ker_source}).
\end{theorem}
The proof of Theorem \ref{ker_L2} is similar to those of Theorem \ref{ker_HK}.

\begin{corollary}\label{col_l2}
Denote  $\theta_{\phi,\zeta}(t)=\frac{t\phi(t)}{\zeta(t)}$ and $\bar{\theta}_{\phi,\zeta}=\frac{\sqrt{t}\phi(t)}{\zeta(t)},$ then with probability at least $1-\delta$ it holds
$$
\|J_{p}^{*}(\beta-J_{p}\tilde{\beta})\|_{\HK}\le C\sqrt{\theta_{\phi,\zeta}^{-1}\left(\frac{1}{\sqrt{N}}+\frac{1}{\sqrt{M}}\right)}\phi\left(\theta_{\phi,\zeta}^{-1}\left(\frac{1}{\sqrt{N}}+\frac{1}{\sqrt{M}}\right)\right)\log^{2}\frac{2}{\delta},$$
for $ \alpha=\alpha_{N,M}:=\theta_{\phi,\zeta}^{-1}\left(\frac{1}{\sqrt{N}}+\frac{1}{\sqrt{M}}\right),$
and 
$$
\|\beta-J_{p}\tilde{\beta}\|_{L_{2,p}}\le C\phi\left(\bar{\theta}_{\phi,\zeta}^{-1}\left(\frac{1}{\sqrt{N}}+\frac{1}{\sqrt{M}}\right)\right)\log^{2}\frac{2}{\delta}, $$
for $\alpha=\alpha_{N,M}:=\bar{\theta}_{\phi,\zeta}^{-1}\left(\frac{1}{\sqrt{N}}+\frac{1}{\sqrt{M}}\right).$

In addition, if $\beta$ meets the source condition (\ref{source_cond:beta2}) with $\phi(t)=t^{s},\, s\in(0,\frac{1}{2} ],$ and $\K$ satisfies Assumption \ref{source_ker} with $\zeta(t)=t^{r},\, r\in(0,\frac{1}{2}],$ then 
\begin{eqnarray}\nonumber
\|J_{p}^{*}(\beta-J_{p}\tilde{\beta})\|_{\HK}\le C\left(\frac{1}{\sqrt{N}}+\frac{1}{\sqrt{M}}\right)^{\frac{1+2s}{2(s-r+1)}}\log^{2}\frac{2}{\delta},
\end{eqnarray}
\end{corollary}
\begin{equation}\nonumber
\|\beta-J_{p}\tilde{\beta}\|_{L_{2,p}}\le C\left(\frac{1}{\sqrt{N}}+\frac{1}{\sqrt{M}}\right)^{\frac{2s}{2(s-r)+1}}\log^{2}\frac{2}{\delta}.
\end{equation}
The proof of Corollary \ref{col_l2}  is similar to those of Corollary \ref{col_HK}.

\section{Computational Cost}
Let's calculate a computational cost which is the number of arithmetic operations required for constructing the approximant $\tilde{\beta}_{M,N,\Znu}^{\alpha_{M,N}}\in\HK^{\Znu}$ within the framework of the method (\ref{eq_beta}).
Note that $\tilde{\beta}_{M,N,\Znu}^{\alpha_{M,N}}\in\HK^{\Znu}$ can be computed with a computational cost $O(|\Znu|^3)$,  which is the computational complexity of solving the system of linear equations (for more details see Appendix C). 

By means of (\ref{sample_choice}):
\begin{equation}\nonumber
|\z^{\nu}|\ge C\mathcal{N}_{\infty}(\alpha)\log\frac{1}{\alpha}\log\frac{1}{\delta},
\end{equation}
we have
\begin{eqnarray}\nonumber
\cost(\tilde{\beta}_{M,N,\Znu}^{\alpha_{M,N}}) = O( N|\Znu|^2 )=O\left(N\left(\mathcal{N}_{\infty}(\alpha)\log\frac{1}{\alpha}\right)^2\right).
\end{eqnarray}
In the scope of the standard assumption that 
\begin{equation}\label{stand:cond}
\mathcal{N}(\alpha)\asymp\alpha^{-s},\quad s\in(0,1],
\end{equation}
 and by Lemma \ref{lem:dem} with $\zeta(t)=t^{\frac{\gamma}{2}},\,\gamma\in(0,1],$ i.e. $\mathcal{N}_{\infty}(\alpha)\le C\alpha^{\gamma-1}$ it follows
 $$\mathcal{N}(\alpha)\le\mathcal{N}_{\infty}(\alpha)\quad \Rightarrow \quad\gamma+s\le1. $$ 
 
Let the parameter $\alpha$ be chosen according with the rule $\mathcal{N}(\alpha)=\alpha\cdot N.$  Then, by (\ref{stand:cond}) for $s\in(0,1-\gamma]$ it holds
$\alpha\asymp N^{-\frac{1}{s+1}}.$ From here, 
\begin{eqnarray}\nonumber
\begin{split}
\cost(\tilde{\beta}_{M,N,\Znu}^{\alpha_{M,N}})&= O\left(N\left(\mathcal{N}_{\infty}(\alpha)\log\frac{1}{\alpha}\right)^2\right)= O\left( N \cdot N^{\frac{2(1-\gamma)}{s+1}}\log^{2}N \right)
= O\left(  N^{\frac{3+s -2\gamma}{s+1}}\log^{2}N \right),
\end{split}
\end{eqnarray}
and hence, the computational cost for computing the Nystr\"om approximant  $\tilde{\beta}_{M,N,\Znu}^{\alpha_{M,N}}\in\HK^{\Znu}$ is subquadratic for $2\gamma+s>1.$ Thus, we proved the following statement. 
\begin{theorem}
Let Assumption \ref{source_ker} is satisfied with $\zeta(t)=t^{\frac{\gamma}{2}},\,\gamma\in(0,1],$ and $\, \mathcal{N}(\alpha)\asymp\alpha^{-s},\quad s\in(0,1-\gamma].$ If $2\gamma+s>1,$ then the Nystr\"om approximant $\tilde{\beta}_{M,N,\Znu}^{\alpha_{M,N}}$ can be computed with subquadratic computational cost.
\end{theorem}

\section{Conclusion}
The present study focuses on the numerical approximation of  the Radon-Nikodym derivative
under the big data assumption.  
To address the above problem, we design an algorithm that is a combination of the Nystr\"om subsampling and the standard  Tikhonov regularization. The convergence rate of the corresponding algorithm is established both in the case when the Radon-Nikodym derivative $\beta$ belongs to RKHS 
$\HK$ and in the case when it does not.  We prove that the proposed approach not only ensure the order of accuracy as algorithms based on the whole sample size, but also allows to achieve subquadratic computational costs in the number of observations.
\\
\\
\text{\bf{\large{Acknowledgments.}}}
We are greatly thankful to the anonymous reviewers for the constructive comments and valuable suggestions, which have significantly improved the clarity of the manuscript.

The second author has received funding through the MSCA4Ukraine project, which is funded by the European Union (ID number 1232599).

\appendix
\newpage

\renewcommand{\thetheorem}{A.\arabic{theorem}}
\renewcommand{\theproposition}{A.\arabic{proposition}}
\renewcommand{\thedefinition}{A.\arabic{definition}}
\renewcommand{\thecorollary}{A.\arabic{corollary}}
\renewcommand{\thelemma}{A.\arabic{lemma}}
\renewcommand{\theremark}{A.\arabic{remark}}
\renewcommand{\theexample}{A.\arabic{example}}
\renewcommand{\theequation}{A.\arabic{equation}}
\section{Appendix. Proof of Auxiliary Statements}\label{ap:a}

{\bf{Proof of Lemma \ref{lem:polar}.}}\\
To prove the first inequality  we apply the polar decomposition of  the operator $J_{p}^{*},$ namely:
$$J_{p}^{*}=U_1^{*}(J_{p}J_{p}^{*})^{\frac{1}{2}},$$
where $U_1\colon \HK\rightarrow L_{2,p},$  $U_1^{*}\colon L_{2,p}\rightarrow\HK,$ which used to call partial isometry operators (see, e.g. \cite[p. 35]{VV}). Moreover $\|U_1\|_{\HK\rightarrow L_{2,p}}=1$, $\|U_1^{*}\|_{ L_{2,p} \rightarrow\HK}=1.$ Hence,
\begin{eqnarray}\nonumber
\begin{split}
\|J_{p}^{*}(\alpha I + J_{p}J_{p}^{*})^{-\frac{1}{2}}\|_{L_{2,p}\rightarrow\HK}&=\|
U_1^{*}(J_{p}J_{p}^{*})^{\frac{1}{2}}(\alpha I + J_{p}J_{p}^{*})^{-\frac{1}{2}}\|_{L_{2,p}\rightarrow\HK}\le\|(J_{p}J_{p}^{*})^{\frac{1}{2}}(\alpha I + J_{p}J_{p}^{*})^{-\frac{1}{2}}\|_{L_{2,p}\rightarrow L_{2,p}}\\
&\le \underset{0<t<l}{\sup}\left(\frac{t}{\alpha+t}\right)^{\frac{1}{2}}\le 1.
\end{split}
\end{eqnarray}
The first inequality is proved. 
Further, using the above reasoning, it is easy to prove the second inequality of the Lemma. Thus, Lemma is proved.\\\\
{\bf{Proof of Lemma \ref{lem:add}.}}\\
Let $f=\phi(J_{p}^{*}J_{p})\mu,$ with $ \phi\in \F, \, \|\mu\|_{\HK}\le 1,$ then
\begin{eqnarray}\label{eq:norm_add}
\begin{split}
\|(\alpha I+J_p^{*}J_p)^{-1/2}f\|_{\HK}&=\|(\alpha I+J_p^{*}J_p)^{-1/2}\phi(J_{p}^{*}J_{p})\mu_{f}\|_{\HK}\le \underset{0<t\leq l}\sup\left((\alpha+t)^{-1}\phi^{2}(t)\right)^{\frac{1}{2}}\\
&\le \frac{1}{\sqrt{\alpha}}\, \underset{0<t\leq l}\sup|\left(1-(\alpha+t)^{-1}t\right)\phi^{2}(t)|^{\frac{1}{2}}.
\end{split}
\end{eqnarray}
Due to the concavity of $\phi^2(t),$ $\phi^2(0)=0,$ for any $t<t_1,$ the point $(t,\phi^{2}(t))\in \mathbb{R}^2$ is above the straight line $v(t)=\frac{\phi^2(t_1)}{t_1}t$ that interpolates the function $\phi^2(t)$ at the points $t=0$ and $t=t_1,$ i.e.
$\phi^2(t)\ge \frac{\phi^2(t_1)}{t_1}t.$ From here we have $\frac{t_1}{\phi^2(t_1)}\ge \frac{t}{\phi^{2}(t)}, \, t_1>t.$ This means, that $\phi^{2}(t)$ is covered by the qualification $p=1.$ Thus, applying Proposition \ref{cover}  to (\ref{eq:norm_add}) we obtain the estimate (\ref{norm_forHK}).
Using the same arguments as for establishing (\ref{norm_forHK}), it is easy  to prove the inequality (\ref{norm_forL2}). Thus, Lemma \ref{lem:add} is proved. \\ \\
{\bf{Proof of Lemma \ref{lem:add0}.}}\\
The proof of the inequality (\ref{eq:prob2}) is given in \cite[Lemmas 2, 8]{RudiComRos15}.
The proof of (\ref{eq:prob20}) is based on (\ref{eq:prob2}) and the obvious equality
\begin{eqnarray}\label{obvious:eq}
\begin{split}
(\alpha I &+ Z)^{\frac{1}{2}} (\alpha I +\PZnu Z\PZnu)^{-1}\PZnu
(\alpha I + Z)^{\frac{1}{2}}
=(\alpha I + Z)^{\frac{1}{2}}\PZnu(\alpha I +\PZnu Z\PZnu)^{-1}\PZnu
(\alpha I + Z)^{\frac{1}{2}}.
\end{split}
\end{eqnarray}
\\ 
{\bf{Proof of Lemma \ref{lem:add01}.}}\\
Denote by
\begin{eqnarray}\label{eq:prob21}
G:=(\alpha I+S_{\X_{p}}^{*}S_{\X_{p}})(\alpha I+\PZnu S_{\X_{p}}^{*}S_{\X_{p}}\PZnu)^{-1}\PZnu.\nonumber
\end{eqnarray}
In view of (\ref{obvious:eq}) we rewrite $G$ as follows
$$
G=(\alpha I+S_{\X_{p}}^{*}S_{\X_{p}})\PZnu(\alpha I+\PZnu S_{\X_{p}}^{*}S_{\X_{p}}\PZnu)^{-1}\PZnu
$$
and give a decomposition
$$
G=G_1+G_2,
$$
where
\begin{eqnarray}\nonumber
\begin{split}
G_1&:=\PZnu(\alpha I+S_{\X_{p}}^{*}S_{\X_{p}})\PZnu(\alpha I+\PZnu S_{\X_{p}}^{*}S_{\X_{p}}\PZnu)^{-1}\PZnu,\\
G_2&:=(I-\PZnu)(\alpha I+S_{\X_{p}}^{*}S_{\X_{p}})\PZnu(\alpha I+\PZnu S_{\X_{p}}^{*}S_{\X_{p}}\PZnu)^{-1}\PZnu.
\end{split}
\end{eqnarray}
It is easy to see that 
$$
G_1=\PZnu(\alpha I+S_{\X_{p}}^{*}S_{\X_{p}})\PZnu(\alpha I+\PZnu S_{\X_{p}}^{*}S_{\X_{p}}\PZnu)^{-1}\PZnu=\PZnu.
$$
Hence, $\|G_1\|_{\HK\rightarrow\HK}=1.$
\\
Now, we are going to bound the norm of $G_2$. To this end, we decompose it  as follows $G_2=G_{21}+G_{22},$  where
\begin{eqnarray}\nonumber
\begin{split}
G_{21}&:=\alpha(I -\PZnu)\PZnu(\alpha I+\PZnu S_{\X_{p}}^{*}S_{\X_{p}}\PZnu)^{-1}\PZnu=0,\\
G_{22}&:=(I-\PZnu)S_{\X_{p}}^{*}S_{\X_{p}}\PZnu(\alpha I+\PZnu S_{\X_{p}}^{*}S_{\X_{p}}\PZnu)^{-1}\PZnu.
\end{split}
\end{eqnarray}
The rest of the proof is about the estimation of the norm of $G_{22}.$
Applying the polar decomposition to the operators $S^{*}_{\X_{p}}$ and $S_{\X_{p}}\PZnu,$ we obtain
\begin{eqnarray}\nonumber
\begin{split}
\|G_{22}\|_{\HK}&\le \|(I-\PZnu)(S^{*}_{\X_{p}}S_{\X_{p}})^{\frac{1}{2}}\|_{\HK\rightarrow\HK}
\|(\PZnu S_{\X_{p}}^{*}S_{\X_{p}}\PZnu)^{\frac{1}{2}}(\alpha I+\PZnu S_{\X_{p}}^{*}S_{\X_{p}}\PZnu)^{-1}\|_{\HK\rightarrow\HK}\\
&\le\|(I-\PZnu)(\alpha I+J^{*}_{p}J_{p})^{\frac{1}{2}}\|_{\HK\rightarrow\HK}\\
&\times \|(\alpha I+J^{*}_{p}J_{p})^{-\frac{1}{2}}(\alpha I+S^{*}_{\X_{p}}S_{\X_{p}})^{\frac{1}{2}}\|_{\HK\rightarrow\HK}\\
&\times\|(\alpha I+S^{*}_{\X_{p}}S_{\X_{p}})^{-\frac{1}{2}}(S^{*}_{\X_{p}}S_{\X_{p}})^{\frac{1}{2}}\|_{\HK\rightarrow\HK}\|(\PZnu S_{\X_{p}}^{*}S_{\X_{p}}\PZnu)^{\frac{1}{2}}(\alpha I+\PZnu S_{\X_{p}}^{*}S_{\X_{p}}\PZnu)^{-1}\|_{\HK\rightarrow\HK}.
\end{split}
\end{eqnarray}
By means of (\ref{eq:g}), (\ref{eq:gp}), (\ref{eq:prob10}), (\ref{bound:3}), and Proposition \ref{Cordes},
we finally get 
$$
\|G_2\|_{\HK}\le C\Big[\frac{B_{N,\alpha}\log\frac{2}{\delta}}{\sqrt{\alpha}}+1\Big]^{\frac{1}{2}}.
$$
Combining obtained estimates, we derive to the statement of Lemma.\\ \\
{\bf{Proof of Lemma \ref{lem:b_L2}.}}\\
Let's consider the following random variable
$$
\xi_{x}(\cdot)=(\alpha I + J_{p}^{*}J_{p})^{-\frac{1}{2}}K(\cdot, x), \quad x\in \X.
$$
It is clear that
$$\|\xi_{x}\|_{\HK}=\|(\alpha I + J_{p}^{*}J_{p})^{-\frac{1}{2}}K(\cdot, x)\|_{\HK}=\sqrt{\mathcal{N}_x(\alpha)}\le \sqrt{\mathcal{N}_{\infty}(\alpha)},$$
$$
\mathbb{E}\xi_{x}=\int_{\X}\, (\alpha I + J_{p}^{*}J_{p})^{-\frac{1}{2}}K(\cdot, x) dq(x)=\int_{\X}\, (\alpha I + J_{p}^{*}J_{p})^{-\frac{1}{2}}K(\cdot, x)\beta(x) dp(x)=(\alpha I + J_{p}^{*}J_{p})^{-\frac{1}{2}}J_{p}^{*}\beta,
$$
and
\begin{eqnarray}\nonumber
\begin{split}
\mathbb{E}\|\xi_{x}\|^{2}_{\HK}&=\int_{\X}\, \|(\alpha I + J_{p}^{*}J_{p})^{-\frac{1}{2}}K(\cdot, x)\|^2_{\HK} dq(x)\le\mathcal{N}_{\infty}(\alpha).
\end{split}
\end{eqnarray}
Thus, for $\xi_
{x^{'}_{j}}(\cdot),$ $j=1,2,\ldots,M,$ drawn i.i.d. from $q(x)$ the conditions of Proposition \ref{concentrat} are satisfied with $L=2\sqrt{\mathcal{N}_{\infty}(\alpha)}$ and $\sigma=\sqrt{\mathcal{N}_{\infty}(\alpha)}.$ Hence, with probability at least $1-\delta$ it holds
\begin{eqnarray}\label{eq_prob}
\begin{split}
\|(\alpha I + J_{p}^{*}J_{p})^{-\frac{1}{2}}( J_{p}^{*}\beta- S_{\X_{q}}^{*}S_{\X_{q}}\textbf{1})\|_{\HK}&=\|\mathbb{E}\xi_{x}-\frac{1}{M}\sum_{j=1}^{M}\,\xi_{x_{j}^{'}}\|_{\HK}
\le C\left(\frac{\sqrt{\mathcal{N}_{\infty}(\alpha)}}{M}+\frac{\sqrt{\mathcal{N}_{\infty}(\alpha)}}{\sqrt{M}}\right)\log^{\frac{1}{2}}\frac{2}{\delta}.\nonumber
\end{split}
\end{eqnarray}
Lemma is proved.

\newpage

\renewcommand{\thetheorem}{B.\arabic{theorem}}
\renewcommand{\theproposition}{B.\arabic{proposition}}
\renewcommand{\thedefinition}{B.\arabic{definition}}
\renewcommand{\thecorollary}{B.\arabic{corollary}}
\renewcommand{\thelemma}{B.\arabic{lemma}}
\renewcommand{\theremark}{B.\arabic{remark}}
\renewcommand{\theexample}{B.\arabic{example}}
\renewcommand{\theequation}{B.\arabic{equation}}

\section{Appendix. Proof of Theorem \ref{L_2}}\label{ap_b}
We split the proof of the theorem into two steps.\\
\text{\bf{Step 1}}  \text{\bf{(Estimation in the metric of $\HK$).}}\\
We start with  the decomposition
\begin{eqnarray}\label{beta_decomp}
\begin{split}
J_{p}^{*}(\beta-J_{p}\tilde{\beta})&=J_{p}^{*}(\beta-J_{p}(\alpha I +\PZnu S_{\X_{p}}^{*}S_{\X_{p}}\PZnu)^{-1}\PZnu S_{\X_{q}}^{*}S_{\X_{q}}\textbf{1})=\overline{\omega}_1+\overline{\omega}_2+\overline{\omega}_3,\nonumber
\end{split}
\end{eqnarray}
where
\begin{eqnarray}\nonumber
\begin{split}
\overline{\omega}_{1}&:=J_p^{*}(I-J_p(\alpha I +\PZnu J_{p}^{*}J_{p}\PZnu)^{-1}\PZnu J_{p}^{*})\beta;\\
\overline{\omega}_2&:=J_{p}^{*}J_{p}\left[(\alpha I +\PZnu J_{p}^{*}J_{p}\PZnu)^{-1}\PZnu J_{p}^{*}\beta-(\alpha I +\PZnu S_{\X_{p}}^{*}S_{\X_{p}}\PZnu)^{-1}\PZnu J_{p}^{*}\beta\right];\\
\overline{\omega}_3&:=J_{p}^{*}J_{p}\left[(\alpha I +\PZnu S_{\X_{p}}^{*}S_{\X_{p}}\PZnu)^{-1}\PZnu J_{p}^{*}\beta-(\alpha I +\PZnu S_{\X_{p}}^{*}S_{\X_{p}}\PZnu)^{-1} \PZnu S_{\X_{q}}^{*}S_{\X_{q}}\textbf{1}\right];\\
\end{split}
\end{eqnarray}

We estimate the norm of each $\overline{\omega}_{i},$  $i=\overline{1,4}.$ For $\overline{\omega}_{1}$ we have
\begin{eqnarray}
\begin{split}
\overline{\omega}_1&=J_{p}^{*}(I-J_{p}(\alpha I +\PZnu J_{p}^{*}J_{p}\PZnu)^{-1}\PZnu J_{p}^{*})\beta
=J_{p}^{*}\left(I-J_{p}\PZnu J_{p}^{*}(\alpha I +J_{p}\PZnu J_{p}^{*})^{-1}\right)\beta\\
&=J_{p}^{*}(\alpha I + J_{p}J_{p}^{*})^{-\frac{1}{2}}
(\alpha I + J_{p}J_{p}^{*})^{\frac{1}{2}}(\alpha I +J_{p}\PZnu J_{p}^{*})^{-\frac{1}{2}}\\
&\times(\alpha I +J_{p}\PZnu J_{p}^{*})^{\frac{1}{2}}\left(I-J_{p}\PZnu J_{p}^{*}(\alpha I +J_{p}\PZnu J_{p}^{*})^{-1}\right)(\alpha I +J_{p}\PZnu J_{p}^{*})^{\frac{1}{2}}\\
&\times(\alpha I +J_{p}\PZnu J_{p}^{*})^{-\frac{1}{2}}(\alpha I + J_{p}J_{p}^{*})^{\frac{1}{2}}(\alpha I + J_{p}J_{p}^{*})^{-\frac{1}{2}}\phi(J_{p}J_{p}^{*})\mu_{\beta}.
\nonumber
\end{split}
\end{eqnarray}
Meanwhile, note that 
\begin{equation}\label{b_gp}
I-J_{p}\PZnu J_{p}^{*}(\alpha I +J_{p}\PZnu J_{p}^{*})^{-1}=\alpha(\alpha I +J_{p}\PZnu J_{p}^{*})^{-1},
\end{equation}
then,
\begin{eqnarray}
\begin{split}
\|\overline{\omega}_1\|_{\HK}&\le\alpha \|J_{p}^{*}(\alpha I + J_{p}J_{p}^{*})^{-\frac{1}{2}}\|_{L_{2,p}\rightarrow\HK}
\|(\alpha I + J_{p}J_{p}^{*})^{\frac{1}{2}}(\alpha I +J_{p}\PZnu J_{p}^{*})^{-\frac{1}{2}}\|_{L_{2,p}\rightarrow L_{2,p} }\\
&\times\|(\alpha I +J_{p}\PZnu J_{p}^{*})^{-\frac{1}{2}}(\alpha I + J_{p}J_{p}^{*})^{\frac{1}{2}}\|_{L_{2,p}\rightarrow L_{2,p} } \|(\alpha I + J_{p}J_{p}^{*})^{-\frac{1}{2}}\phi(J_{p}J_{p}^{*})\mu_{\beta}\|_{L_{2,p} }.
\nonumber
\end{split}
\end{eqnarray}
Moreover,
\begin{eqnarray}\nonumber
\begin{split}
\|(\alpha I &+J_{p}\PZnu J_{p}^{*})^{-1}(\alpha I + J_{p}J_{p}^{*})\|_{L_{2,p}\rightarrow L_{2,p} }\le \|(\alpha I +J_{p}\PZnu J_{p}^{*})^{-1}\left(J_{p}J_{p}^{*}-J_{p}\PZnu J_{p}^{*}\right)\|_{L_{2,p}\rightarrow L_{2,p} }+1\\
&\le \|(\alpha I +J_{p}\PZnu J_{p}^{*})^{-1}\|_{L_{2,p}\rightarrow L_{2,p}} 
\|J_p(I-\PZnu)\|^{2}_{\HK\rightarrow L_{2,p}}+1.
\end{split}
\end{eqnarray}
By means of (\ref{eq:g}) and (\ref{eq:prob1}) we get
\begin{eqnarray}\label{b_add_err}
\begin{split}
\|(\alpha I& +J_{p}\PZnu J_{p}^{*})^{-1}(\alpha I + J_{p}J_{p}^{*})\|_{L_{2,p}\rightarrow L_{2,p}} \le \frac{\tilde{\gamma}}{\alpha} \cdot 3\alpha+1=3\tilde{\gamma}+1.
\end{split}
\end{eqnarray}
This together with  Lemma \ref{lem:polar}, Proposition \ref{Cordes}, and (\ref{norm_forL2}) implies that
\begin{eqnarray}\label{om:1}
\begin{split}
\|\overline{\omega}_1\|_{\HK}&\le 
C\sqrt{\alpha}\phi(\alpha).
\end{split}
\end{eqnarray}

Now, we are going to bound the norm of $\overline{\omega}_{2}:$ 
\begin{eqnarray}\nonumber
\begin{split}
\overline{\omega}_2&=J_{p}^{*}J_{p}\left[(\alpha I +\PZnu J_{p}^{*}J_{p}\PZnu)^{-1}-(\alpha I +\PZnu S_{\X_{p}}^{*}S_{\X_{p}}\PZnu)^{-1}\right]\PZnu J_{p}^{*}\beta\\
&=J_{p}^{*}J_{p}(\alpha I +\PZnu S_{\X_{p}}^{*}S_{\X_{p}}\PZnu)^{-1}\PZnu\left[ S_{\X_{p}}^{*}S_{\X_{p}}-J_{p}^{*}J_{p}\right]\PZnu (\alpha I +\PZnu J_{p}^{*}J_{p}\PZnu)^{-1}\PZnu J_{p}^{*}\beta\\
&=J_{p}^{*}J_{p}(\alpha I+J_{p}^{*}J_{p})^{-1}(\alpha I+J_{p}^{*}J_{p})(\alpha I + S_{\X_{p}}^{*}S_{\X_{p}})^{-1}\\
&\times(\alpha I + S_{\X_{p}}^{*}S_{\X_{p}})
(\alpha I +\PZnu S_{\X_{p}}^{*}S_{\X_{p}}\PZnu)^{-1}\PZnu\\
&\times\left[ S_{\X_{p}}^{*}S_{\X_{p}}-J_{p}^{*}J_{p}\right](\alpha I + J_{p}^{*}J_{p})^{-\frac{1}{2}}\\
&\times(\alpha I + J_{p}^{*}J_{p})^{\frac{1}{2}}\PZnu(\alpha I +\PZnu J_{p}^{*}J_{p}\PZnu)^{-1}\PZnu (\alpha I + J_{p}^{*}J_{p})^{\frac{1}{2}}(\alpha I + J_{p}^{*}J_{p})^{-\frac{1}{2}}J_{p}^{*}\beta,\\
\end{split}
\end{eqnarray}
then
\begin{eqnarray}\nonumber
\begin{split}
\|\overline{\omega}_2\|_{\HK}
&\le
\|J_p^{*}J_p(\alpha I+J_{p}^{*}J_{p})^{-1}\|_{\HK \rightarrow \HK}\\
&\times\|(\alpha I+J_{p}^{*}J_{p})(\alpha I + S_{\X_{p}}^{*}S_{\X_{p}})^{-1}\|_{\HK \rightarrow \HK}\\
&\times\|(\alpha I + S_{\X_{p}}^{*}S_{\X_{p}})
(\alpha I +\PZnu S_{\X_{p}}^{*}S_{\X_{p}}\PZnu)^{-1}\PZnu\|_{\HK \rightarrow \HK}\\
&\times\|\left[ S_{\X_{p}}^{*}S_{\X_{p}}-J_{p}^{*}J_{p}\right](\alpha I + J_{p}^{*}J_{p})^{-\frac{1}{2}}\|_{\HK \rightarrow \HK}\\
&\times\|(\alpha I + J_{p}^{*}J_{p})^{\frac{1}{2}}\PZnu(\alpha I +\PZnu J_{p}^{*}J_{p}\PZnu)^{-1}\PZnu (\alpha I + J_{p}^{*}J_{p})^{\frac{1}{2}}\|_{\HK \rightarrow \HK}\\
&\times\|J_{p}^{*}(\alpha I + J_{p}J_{p}^{*})^{-\frac{1}{2}}\|_{L_{2,p} \rightarrow \HK}\|\beta\|_{L_{2,p}}.\\
\end{split}
\end{eqnarray}
By (\ref{eq:gp}), (\ref{eq:prob2}) with $Z=J_{p}^{*}J_{p}$, (\ref{bound:1}), (\ref{bound:2}), Lemma \ref{lem:add01}, Proposition \ref{Cordes}, and Lemma \ref{lem:polar}, we get
\begin{eqnarray}\nonumber
\begin{split}
\|\overline{\omega}_2\|_{\HK}
&\le C\left(1+C\left(\frac{\mathcal{B}_{N,\alpha}\log\frac{2}{\delta}}{\sqrt{\alpha}}+1\right)^{\frac{1}{2}}\right)\left(\left(\frac{B_{N,\alpha}\log\frac{2}{\delta}}{\sqrt{\alpha}}\right)^2+1\right)
B_{N,\alpha}\log\frac{2}{\delta}.
\end{split}
\end{eqnarray}
In view of (\ref{est:B}), we have
\begin{eqnarray}\label{om:2}
\begin{split}
\|\overline{\omega}_2\|_{\HK}\le C\left(1+C\left(\frac{\mathcal{B}_{N,\alpha}\log\frac{2}{\delta}}{\sqrt{\alpha}}+1\right)^{\frac{1}{2}}\right)\left(\left(\frac{B_{N,\alpha}\log{\frac{2}{\delta}}}{\sqrt{\alpha}}\right)^2+1\right)\frac{\sqrt{\mathcal{N}_{\infty}(\alpha)}}{\sqrt{N}}\log\frac{2}{\delta}.
\end{split}
\end{eqnarray}

We are at the point to bound the norm of $\overline{\omega}_{3}.$ We start with the decomposition
\begin{eqnarray}\nonumber
\begin{split}
\overline{\omega}_3&=J_{p}^{*}J_{p}(\alpha I +\PZnu S_{\X_{p}}^{*}S_{\X_{p}}\PZnu)^{-1}\PZnu( J_{p}^{*}\beta- S_{\X_{q}}^{*}S_{\X_{q}}\textbf{1})\\
&=J_{p}^{*}J_{p}(\alpha I+J_{p}^{*}J_{p})^{-\frac{1}{2}}(\alpha I+J_{p}^{*}J_{p})^{\frac{1}{2}}(\alpha I + S_{\X_{p}}^{*}S_{\X_{p}})^{-\frac{1}{2}}\\
&\times(\alpha I + S_{\X_{p}}^{*}S_{\X_{p}})^{\frac{1}{2}}
(\alpha I +\PZnu S_{\X_{p}}^{*}S_{\X_{p}}\PZnu)^{-1}\PZnu
(\alpha I + S_{\X_{p}}^{*}S_{\X_{p}})^{\frac{1}{2}}\\
&\times(\alpha I + S_{\X_{p}}^{*}S_{\X_{p}})^{-\frac{1}{2}}(\alpha I + J_{p}^{*}J_{p})^{\frac{1}{2}}(\alpha I + J_{p}^{*}J_{p})^{-\frac{1}{2}}( J_{p}^{*}\beta- S_{\X_{q}}^{*}S_{\X_{q}}\textbf{1}),
\end{split}
\end{eqnarray}
then
\begin{eqnarray}\nonumber
\begin{split}
\|\overline{\omega}_3\|_{\HK}
&\le
\|(J_p^{*}J_p)^{\frac{1}{2}}\|_{\HK\rightarrow\HK}
\|(J_p^{*}J_p)^{\frac{1}{2}}(\alpha I+J_{p}^{*}J_{p})^{-\frac{1}{2}}\|_{\HK \rightarrow \HK}
\|(\alpha I+J_{p}^{*}J_{p})^{\frac{1}{2}}(\alpha I + S_{\X_{p}}^{*}S_{\X_{p}})^{-\frac{1}{2}}\|_{\HK \rightarrow \HK}\\
&\times\|(\alpha I + S_{\X_{p}}^{*}S_{\X_{p}})^{\frac{1}{2}}
(\alpha I +\PZnu S_{\X_{p}}^{*}S_{\X_{p}}\PZnu)^{-1}\PZnu
(\alpha I + S_{\X_{p}}^{*}S_{\X_{p}})^{\frac{1}{2}}\|_{\HK \rightarrow \HK}\\
&\times\|(\alpha I + S_{\X_{p}}^{*}S_{\X_{p}})^{-\frac{1}{2}}(\alpha I + J_{p}^{*}J_{p})^{\frac{1}{2}}\|_{\HK \rightarrow \HK}\|(\alpha I + J_{p}^{*}J_{p})^{-\frac{1}{2}}( J_{p}^{*}\beta- S_{\X_{q}}^{*}S_{\X_{q}}\textbf{1})\|_{\HK}.
\end{split}
\end{eqnarray}
Applying Lemma \ref{lem:polar}, (\ref{eq:prob20}) with $Z=S_{\X_{p}}^{*}S_{\X_{p}}$, (\ref{bound:2}), Proposition \ref{Cordes}, and  Lemma \ref{lem:b_L2}, we derive
\begin{eqnarray}\label{om:3}
\begin{split}
\|\overline{\omega}_3\|_{\HK}
&\le C\left(\left(\frac{B_{N,\alpha}\log{\frac{2}{\delta}}}{\sqrt{\alpha}}\right)^2+1\right)\frac{\sqrt{\mathcal{N}_{\infty}(\alpha)}}{\sqrt{M}}\log\frac{2}{\delta}.
\end{split}
\end{eqnarray}
Summing up (\ref{om:1}), (\ref{om:2}), and (\ref{om:3}), with probability at least $1-\delta$ we finally get
 \begin{eqnarray}\nonumber
 \begin{split}
\|J_{p}^{*}(\beta&-J_{p}\tilde{\beta})\|_{\HK}\le C\sqrt{\alpha}\phi(\alpha)
\\
&+ 
C\left(1+C\left(\frac{\mathcal{B}_{N,\alpha}\log\frac{2}{\delta}}{\sqrt{\alpha}}+1\right)^{\frac{1}{2}}\right)\left[\left(\frac{B_{N,\alpha}\log\frac{2}{\delta}}{\sqrt{\alpha}}\right)^{2}+1\right]\left(\frac{1}{\sqrt{N}}+
\frac{1}{\sqrt{M}}\right)\sqrt{\mathcal{N}_{\infty}(\alpha)}\log\frac{2}{\delta}.
\end{split}
\end{eqnarray}
\\
\\
\text{\bf{Step 2}} \text{\bf{(Estimation in the metric of $L_{2,p}$).}}\\
Similarly to \text{\bf{Step 1}}, we start with the decomposition
\begin{eqnarray}\label{beta_decomp}
\begin{split}
\beta-J_{p}\tilde{\beta}&=\beta-J_{p}(\alpha I +\PZnu S_{\X_{p}}^{*}S_{\X_{p}}\PZnu)^{-1}\PZnu S_{\X_{q}}^{*}S_{\X_{q}}\textbf{1}=\overline{\sigma}_1+\overline{\sigma}_2+\overline{\sigma}_3,\nonumber
\end{split}
\end{eqnarray}
where
\begin{eqnarray}\nonumber
\begin{split}
\overline{\sigma}_{1}&:=(I-J_p(\alpha I +\PZnu J_{p}^{*}J_{p}\PZnu)^{-1}\PZnu J_{p}^{*})\beta;\\
\overline{\sigma}_2&:=J_{p}\left[(\alpha I +\PZnu J_{p}^{*}J_{p}\PZnu)^{-1}\PZnu J_{p}^{*}\beta-(\alpha I +\PZnu S_{\X_{p}}^{*}S_{\X_{p}}\PZnu)^{-1}\PZnu J_{p}^{*}\beta\right];\\
\overline{\sigma}_3&:=J_{p}\left[(\alpha I +\PZnu S_{\X_{p}}^{*}S_{\X_{p}}\PZnu)^{-1}\PZnu J_{p}^{*}\beta-(\alpha I +\PZnu S_{\X_{p}}^{*}S_{\X_{P}}\PZnu)^{-1}\PZnu S_{\X_{q}}^{*}S_{\X_{q}}\textbf{1}\right].\\
\end{split}
\end{eqnarray}

We estimate the norm of each $\overline{\sigma}_{i},$  $i=\overline{1,3}.$ Using (\ref{b_gp}), for $\overline{\sigma}_1$ we have
\begin{eqnarray}
\begin{split}
\overline{\sigma}_1&=(I-J_{p}(\alpha I +\PZnu J_{p}^{*}J_{p}\PZnu)^{-1}\PZnu J_{p}^{*})\beta
=\left(I-J_{p}\PZnu J_{p}^{*}(\alpha I +J_{p}\PZnu J_{p}^{*})^{-1}\right)\beta\\
&=\alpha(\alpha I + J_{p}J_{p}^{*})^{-\frac{1}{2}}
(\alpha I + J_{p}J_{p}^{*})^{\frac{1}{2}}(\alpha I +J_{p}\PZnu J_{p}^{*})^{-1}
(\alpha I + J_{p}J_{p}^{*})^{\frac{1}{2}}(\alpha I + J_{p}J_{p}^{*})^{-\frac{1}{2}}\phi(J_{p}J_{p}^{*})\mu_{\beta}.\nonumber
\end{split}
\end{eqnarray}
Applying (\ref{eq:g}), (\ref{b_add_err}), Proposition \ref{Cordes}, and (\ref{norm_forL2}), we get
\begin{eqnarray}\label{sigg:1}
\begin{split}
\|\overline{\sigma}_1\|_{L_{2,p}}&\le\alpha\|(\alpha I + J_{p}J_{p}^{*})^{-\frac{1}{2}}\|_{L_{2,p}\rightarrow L_{2,p}}
\|(\alpha I + J_{p}J_{p}^{*})^{\frac{1}{2}}(\alpha I +J_{p}\PZnu J_{p}^{*})^{-\frac{1}{2}}\|_{L_{2,p}\rightarrow L_{2,p} }\\
&\times\|(\alpha I +J_{p}\PZnu J_{p}^{*})^{-\frac{1}{2}}(\alpha I + J_{p}J_{p}^{*})^{\frac{1}{2}}\|_{L_{2,p}\rightarrow L_{2,p} }
\|(\alpha I + J_{p}J_{p}^{*})^{-\frac{1}{2}}\phi(J_{p}J_{p}^{*})\mu_{\beta}\|_{L_{2,p} }\le C\phi(\alpha).
\end{split}
\end{eqnarray}

Now, we are going to bound the norm of $\overline{\sigma}_{2}.$ 
\begin{eqnarray}\nonumber
\begin{split}
\overline{\sigma}_2&=J_{p}\left[(\alpha I +\PZnu J_{p}^{*}J_{p}\PZnu)^{-1}-(\alpha I +\PZnu S_{\X_{p}}^{*}S_{\X_{p}}\PZnu)^{-1}\right]\PZnu J_{p}^{*}\beta\\
&=J_{p}(\alpha I +\PZnu S_{\X_{p}}^{*}S_{\X_{p}}\PZnu)^{-1}\PZnu\left[ S_{\X_{p}}^{*}S_{\X_{p}}-J_{p}^{*}J_{p}\right]\PZnu (\alpha I +\PZnu J_{p}^{*}J_{p}\PZnu)^{-1}\PZnu J_{p}^{*}\beta\\
&=J_{p}(\alpha I+J_{p}^{*}J_{p})^{-\frac{1}{2}}(\alpha I+J_{p}^{*}J_{p})^{\frac{1}{2}}(\alpha I + S_{\X_{p}}^{*}S_{\X_{p}})^{-\frac{1}{2}}\\
&\times(\alpha I + S_{\X_{p}}^{*}S_{\X_{p}})^{\frac{1}{2}}
(\alpha I +\PZnu S_{\X_{p}}^{*}S_{\X_{p}}\PZnu)^{-1}\PZnu
(\alpha I + S_{\X_{p}}^{*}S_{\X_{p}})^{\frac{1}{2}}\\
&\times(\alpha I + S_{\X_{p}}^{*}S_{\X_{p}})^{-\frac{1}{2}}(\alpha I + J_{p}^{*}J_{p})^{\frac{1}{2}}(\alpha I + J_{p}^{*}J_{p})^{-\frac{1}{2}}\left[ S_{\X_{p}}^{*}S_{\X_{p}}-J_{p}^{*}J_{p}\right]\\
&\times\PZnu(\alpha I +\PZnu J_{p}^{*}J_{p}\PZnu)^{-1}\PZnu J_{p}^{*}\beta,\\
\end{split}
\end{eqnarray}
then
\begin{eqnarray}\nonumber
\begin{split}
\|\overline{\sigma}_2\|_{L_{2,p}}
&\le
\|J_{p}(\alpha I+J_{p}^{*}J_{p})^{-\frac{1}{2}}\|_{\HK \rightarrow L_{2,p}}\\
&\times\|(\alpha I+J_{p}^{*}J_{p})^{\frac{1}{2}}(\alpha I + S_{\X_{p}}^{*}S_{\X_{p}})^{-\frac{1}{2}}\|_{\HK \rightarrow \HK}\\
&\times\|(\alpha I + S_{\X_{p}}^{*}S_{\X_{p}})^{\frac{1}{2}}
(\alpha I +\PZnu S_{\X_{p}}^{*}S_{\X_{p}}\PZnu)^{-1}\PZnu
(\alpha I + S_{\X_{p}}^{*}S_{\X_{p}})^{\frac{1}{2}}\|_{\HK \rightarrow \HK}\\
&\times\|(\alpha I + S_{\X_{p}}^{*}S_{\X_{p}})^{-\frac{1}{2}}(\alpha I + J_{p}^{*}J_{p})^{\frac{1}{2}}\|_{\HK \rightarrow \HK}\|(\alpha I + J_{p}^{*}J_{p})^{-\frac{1}{2}}\left[ S_{\X_{p}}^{*}S_{\X_{p}}-J_{p}^{*}J_{p}\right]\|_{\HK \rightarrow \HK}\\
&\times\|(\alpha I +\PZnu J_{p}^{*}J_{p}\PZnu)^{-1}\PZnu J_{p}^{*}\phi(J_{p}J_{p}^{*})\mu_{\beta}\|_{\HK}.\\
\end{split}
\end{eqnarray}
For the last norm, we have
\begin{eqnarray}\nonumber
\begin{split}
&\|(\alpha I +\PZnu J_{p}^{*}J_{p}\PZnu)^{-1}\PZnu J_{p}^{*}\phi(J_{p}J_{p}^{*})\mu_{\beta}\|_{\HK}=\|\PZnu J_{p}^{*}(\alpha I +J_{p}\PZnu J^{*}_{p})^{-1}\phi(J_{p}J_{p}^{*})\mu_{\beta}\|_{\HK}\\
&\le\|\PZnu J_{p}^{*}(\alpha I +J_{p}\PZnu J^{*}_{p})^{-1}\phi(J_{p}\PZnu J_{p}^{*})\mu_{\beta}\|_{\HK}
+\|\PZnu J_{p}^{*}(\alpha I +J_{p}\PZnu J^{*}_{p})^{-1}\left(\phi( J_{p}J_{p}^{*})-\phi(J_{p}\PZnu J_{p}^{*})\right)\mu_{\beta}\|_{\HK}.
\end{split}
\end{eqnarray}
Keeping in mind that $\phi$ is operator monotone function, by means of 
(\ref{eq:g}), (\ref{qualific_root}), (\ref{eq:prob1}), and (\ref{eq:qual}), we obtain
\begin{eqnarray}\label{est_add}\nonumber
\begin{split}
\|(\alpha I &+\PZnu J_{p}^{*}J_{p}\PZnu)^{-1}\PZnu J_{p}^{*}\phi(J_{p}J_{p}^{*})\mu_{\beta}\|_{\HK}\le \frac{1}{\sqrt{\alpha}}\phi(\alpha)+\frac{\overline{\gamma}}{\sqrt{\alpha}}\phi\left(\|J_{p}(I-\PZnu)\|^{2}_{\HK\rightarrow L_{2,p}}\right)\le\frac{C}{\sqrt{\alpha}}\phi(\alpha).
\end{split}
\end{eqnarray}
Then from this and 
 Lemma \ref{lem:polar}, (\ref{eq:prob20}) with $Z=S_{\X_{p}}^{*}S_{\X_{p}}$, (\ref{bound:1}), (\ref{bound:2}), Proposition \ref{Cordes}, we derive
\begin{eqnarray}\label{sigg:2}
\begin{split}
\|\overline{\sigma}_2\|_{{L_{2,p}}}\le C\left(\left(\frac{B_{N,\alpha}\log{\frac{2}{\delta}}}{\sqrt{\alpha}}\right)^2+1\right)\frac{B_{N,\alpha}}{\sqrt{\alpha}}\log{\frac{2}{\delta}}\phi(\alpha).
\end{split}
\end{eqnarray}

We are at the point to bound the norm of $\overline{\sigma}_{3}.$ We start with the decomposition
\begin{eqnarray}\nonumber
\begin{split}
\overline{\sigma}_3&=J_{p}(\alpha I +\PZnu S_{\X_{p}}^{*}S_{\X_{p}}\PZnu)^{-1}\PZnu( J_{p}^{*}\beta- S_{\X_{q}}^{*}S_{\X_{q}}\textbf{1})\\
&=J_{p}(\alpha I+J_{p}^{*}J_{p})^{-\frac{1}{2}}(\alpha I+J_{p}^{*}J_{p})^{\frac{1}{2}}(\alpha I + S_{\X_{p}}^{*}S_{\X_{p}})^{-\frac{1}{2}}\\
&\times(\alpha I + S_{\X_{p}}^{*}S_{\X_{p}})^{\frac{1}{2}}
(\alpha I +\PZnu S_{\X_{p}}^{*}S_{\X_{p}}\PZnu)^{-1}\PZnu
(\alpha I + S_{\X_{p}}^{*}S_{\X_{p}})^{\frac{1}{2}}\\
&\times(\alpha I + S_{\X_{p}}^{*}S_{\X_{p}})^{-\frac{1}{2}}(\alpha I + J_{p}^{*}J_{p})^{\frac{1}{2}}(\alpha I + J_{p}^{*}J_{p})^{-\frac{1}{2}}( J_{p}^{*}\beta- S_{\X_{q}}^{*}S_{\X_{q}}\textbf{1}),
\end{split}
\end{eqnarray}
then
\begin{eqnarray}\nonumber
\begin{split}
\|\overline{\sigma}_3\|_{L_{2,p}}
&\le
\|J_{p}(\alpha I+J_{p}^{*}J_{p})^{-\frac{1}{2}}\|_{\HK \rightarrow L_{2,p}}
\|(\alpha I+J_{p}^{*}J_{p})^{\frac{1}{2}}(\alpha I + S_{\X_{p}}^{*}S_{\X_{p}})^{-\frac{1}{2}}\|_{\HK \rightarrow \HK}\\
&\times\|(\alpha I + S_{\X_{p}}^{*}S_{\X_{p}})^{\frac{1}{2}}
(\alpha I +\PZnu S_{\X_{p}}^{*}S_{\X_{p}}\PZnu)^{-1}\PZnu
(\alpha I + S_{\X_{p}}^{*}S_{\X_{p}})^{\frac{1}{2}}\|_{\HK \rightarrow \HK}\\
&\times\|(\alpha I + S_{\X_{p}}^{*}S_{\X_{p}})^{-\frac{1}{2}}(\alpha I + J_{p}^{*}J_{p})^{\frac{1}{2}}\|_{\HK \rightarrow \HK}\|(\alpha I + J_{p}^{*}J_{p})^{-\frac{1}{2}}( J_{p}^{*}\beta- S_{\X_{q}}^{*}S_{\X_{q}}\textbf{1})\|_{\HK}.
\end{split}
\end{eqnarray}
Using  (\ref{eq:prob20}) with $Z=S_{\X_{p}}^{*}S_{\X_{p}}$, (\ref{bound:2}), Proposition \ref{Cordes}, Lemma \ref{lem:polar}, and Lemma \ref{lem:b_L2}, we have
\begin{eqnarray}\label{sigg:3}
\begin{split}
\|\overline{\sigma}_3\|_{{L_{2,p}}}
&\le C\left(\left(\frac{B_{N,\alpha}\log{\frac{2}{\delta}}}{\sqrt{\alpha}}\right)^2+1\right)\frac{\sqrt{\mathcal{N}_{\infty}(\alpha)}}{\sqrt{M}}\log^{\frac{1}{2}}\frac{2}{\delta}.
\end{split}
\end{eqnarray}
Summing up (\ref{sigg:1}), (\ref{sigg:2}), and (\ref{sigg:3}), with probability at least $1-\delta$ we finally get
 \begin{eqnarray}\nonumber
 \begin{split}
\|\beta-J_{p}\tilde{\beta}\|_{L_{2,p}}&\le C\phi(\alpha)
+ C\left[\left(\frac{B_{N,\alpha}\log{\frac{2}{\delta}}}{\sqrt{\alpha}}\right)^2+1\right]\frac{B_{N,\alpha}}{\sqrt{\alpha}}\phi(\alpha)\log{\frac{2}{\delta}}\\
&+ C\left[\left(\frac{B_{N,\alpha}\log\frac{2}{\delta}}{\sqrt{\alpha}}\right)^{2}+1\right]
\frac{\sqrt{\mathcal{N}_{\infty}(\alpha)}}{\sqrt{M}}\log^{\frac{1}{2}}\frac{2}{\delta}.
\end{split}
\end{eqnarray}
Thus, Theorem \ref{L_2} is completely proved. 
\newpage

\renewcommand{\thetheorem}{C.\arabic{theorem}}
\renewcommand{\theproposition}{C.\arabic{proposition}}
\renewcommand{\thedefinition}{C.\arabic{definition}}
\renewcommand{\thecorollary}{C.\arabic{corollary}}
\renewcommand{\thelemma}{C.\arabic{lemma}}
\renewcommand{\theremark}{C.\arabic{remark}}
\renewcommand{\theexample}{C.\arabic{example}}
\renewcommand{\theequation}{C.\arabic{equation}}
\section{Appendix. Explanation to (\ref{eq_beta})}
Following \cite{NgPer23}, we present $\tilde{\beta}_{\X}^{\alpha}$ as
\begin{eqnarray}\label{add:beta_Per}
\tilde{\beta}_{\X}^{\alpha}=(\alpha I + S_{\X_{p}}^{*}S_{\X_{p}})^{-1}S_{\X_{q}}^{*}S_{\X_{q}} \textbf{1}.
\end{eqnarray}
Let us rewrite (\ref{add:beta_Per}) as the following equation
\begin{eqnarray}\label{add:beta_eq}
(\alpha I + S_{\X_{p}}^{*}S_{\X_{p}})\tilde{\beta}=S_{\X_{q}}^{*}S_{\X_{q}} \textbf{1},
\end{eqnarray}
where $\tilde{\beta}=\tilde{\beta}_{\X}^{\alpha}.$
Further, based on the representation of the operators $S_{\X_{p}}^{*},\, S_{\X_{q}}^{*}$, a solution of the equation (\ref{add:beta_eq}) we will seek as 
\begin{eqnarray}\nonumber
\tilde{\beta}=\sum_{i=1}^{N}\,c_i\K(\cdot,x_i)+\sum_{j=1}^{M}\,c^{'}_j\K(\cdot,x^{'}_j).
\end{eqnarray}
Thus, substituting $\tilde{\beta}$ into (\ref{add:beta_eq}), we obtain
\begin{eqnarray}\label{SLAE_Per}
\begin{split}
&\alpha \sum_{i=1}^{N}\,c_i\K(\cdot,x_i) +\alpha\sum_{i=1}^{M}\,c^{'}_j\K(\cdot,x^{'}_j)\\
&+\frac{1}{N}\sum_{k=1}^{N}\K(\cdot,x_k) \,\sum_{i=1}^{N}\,c_i\K(x_k,x_i)+\frac{1}{N}\sum_{k=1}^{N}\K(\cdot,x_k) \,\sum_{j=1}^{M}\,c^{'}_j\K(x_k,x^{'}_j)=\frac{1}{M}\sum_{j=1}^{M}\, \K(\cdot,x^{'}_j).
\end{split}
\end{eqnarray}
We assume, that $\K(\cdot,x_i) $ and $\K(\cdot,x^{'}_j) $ are linearly independent. Then, (\ref{SLAE_Per}) derives to  two systems of  linear equations:\\
1) for any $k=1,\ldots, N:$ 
$$\alpha c_{i}+\frac{1}{N}\left(\sum_{i=1}^{N}\,c_i\K(x_k,x_i)+\sum_{j=1}^{M}\,c^{'}_j\K(x_k,x^{'}_j)\right)=0,$$\\
2) for any $j=1,\ldots, M:$
$$
\alpha c_{j}^{'}=\frac{1}{M}.
$$
It follows that $ c_{j}^{'}=\frac{1}{\alpha M}.$ Substituting obtained $c_{j}^{'} $ $\text{into}$ the first system, for all $k=1,\ldots, N$ we get
\begin{eqnarray}\label{SLAE_Per_fin}
\begin{split}
\alpha c_{i}+\frac{1}{N}\sum_{i=1}^{N}\,c_i\K(x_k,x_i)=-\frac{1}{\alpha M N}\sum_{j=1}^{M}\,\K(x_k,x^{'}_j).
\end{split}
\end{eqnarray}

Analysis of the system (\ref{SLAE_Per_fin}) shows that we operate with  two Gram's matrices $\K(x_k,x_i)$ and $\K(x_k,x_j^{'})$ with dimensions $N\times N$,  $N\times M$, correspondingly. From here, it follows that within the framework of the Nystr\"om method samples from $\{x_i\}_{i=1}^{N} $ and $\{x_j^{'}\}_{j=1}^{M}$ can be selected independently, but to reduce the computational cost, samples size should be chosen equal. 

It worth to note that the computational cost of the algorithm from  \cite{NgPer23}, more accurately the cost associated with the computation of the minimizer $\tilde{\beta}$, is $O(N^{3}),$ which is the computational complexity of solving (\ref{SLAE_Per_fin}).
Therefore, in the setting, where $N$ is large enough, it is necessary to avoid the computation of the minimizers $\tilde{\beta}.$

Now, we are going to calculate the computational cost of the minimizer (\ref{eq_beta}).
Recall (see (\ref{HK})) that
\begin{eqnarray}\nonumber
\HK^{\Znu}\colon=\left\{f\colon f=\sum_{i=1}^{|\Znu|}\,c_i\K(\cdot,x_i)+\sum_{j=1}^{|\Znu|}\,c^{'}_j\K(\cdot,x^{'}_j)\right\}, 
\end{eqnarray}
where $|\Znu|\ll \min\{N,M\}.$
In the scope of the regularized Nystr\"om subsampling the approximation to the Radon-Nikodym derivative has the form (see (\ref{eq_beta})):
\begin{eqnarray}\label{add:beta}
\tilde{\beta}_{M,N,,\Znu}^{\alpha_{M,N}}=(\alpha I +\PZnu S_{\X_{p}}^{*}S_{\X_{p}}\PZnu)^{-1}\PZnu S_{\X_{q}}^{*}S_{\X_{q}}\textbf{1}.
\end{eqnarray}
As before, we rewrite  (\ref{add:beta}) as follows
\begin{eqnarray}\label{add:eq}
(\alpha I +\PZnu S_{\X_{p}}^{*}S_{\X_{p}}\PZnu)\tilde{\beta}=\PZnu S_{\X_{q}}^{*}S_{\X_{q}} \textbf{1},
\end{eqnarray}
where $\tilde{\beta}=\tilde{\beta}_{M,N,\Znu}^{\alpha_{M,N}}.$
Based on the representation of the operators $S_{\X_{p}}^{*},\, S_{\X_{q}}^{*}$, a solution $\tilde{\beta}$ we will seek as 
\begin{eqnarray}\nonumber
\tilde{\beta}=\sum_{i=1}^{N}\,c_i\K(\cdot,x_i)+\sum_{j=1}^{M}\,c^{'}_j\K(\cdot,x^{'}_j).
\end{eqnarray}
Further, substituting $\tilde{\beta}$ into (\ref{add:eq}), we derive
\begin{eqnarray}\label{SLAE}
\begin{split}
\alpha \sum_{i=1}^{|\Znu|}\,c_i\K(\cdot,x_i) &+\alpha\sum_{j=1}^{|\Znu|}\,c^{'}_j\K(\cdot,x^{'}_j)\\
&+\frac{1}{N}\sum_{k=1}^{|\Znu|}\K(\cdot,x_k) \,\sum_{i=1}^{|\Znu|}\,c_i\K(x_k,x_i)+\frac{1}{N}\sum_{k=1}^{|\Znu|}\K(\cdot,x_k) \,\sum_{j=1}^{|\Znu|}\,c^{'}_j\K(x_k,x^{'}_j)=\frac{1}{M}\sum_{j=1}^{|\Znu|}\,\K(\cdot,x^{'}_j).
\end{split}
\end{eqnarray}
Assuming that $\K(\cdot,x_i) $ and $\K(\cdot,x^{'}_j) $ are linearly independent,  (\ref{SLAE}) leads to the systems of  linear equations:\\
1) for any $k=1,\ldots, |\Znu|:$ 
$$\alpha c_{i}+\frac{1}{N}\left(\sum_{i=1}^{|\Znu|}\,c_i\K(x_k,x_i)+\sum_{j=1}^{|\Znu|}\,c^{'}_j\K(x_k,x^{'}_j)\right)=0,$$\\
2) for any $j=1,\ldots, |\Znu|\colon$
$
\alpha c_{j}^{'}=\frac{1}{M}\quad \Rightarrow \quad c_{j}^{'}=\frac{1}{\alpha M}.
$
\\
Substituting obtained $c_{j}^{'} $ $\text{into}$ the first system, for all $k=1,\ldots, |\Znu|$, we get
\begin{eqnarray}\label{SLAE_fin}
\begin{split}
\alpha c_{i}+\frac{1}{N}\sum_{i=1}^{|\Znu|}\,c_i\K(x_k,x_i)=-\frac{1}{\alpha M N}\sum_{j=1}^{|\Znu|}\,\K(x_k,x^{'}_j).
\end{split}
\end{eqnarray}
From (\ref{SLAE_fin}) it follows that the computational cost that is needed to design  the Nystr\"om approximant  (\ref{add:beta}) is of order $O(|\Znu|^{3}),$ which in turn is much less than $O(N^{3}).$ Thus, proposed approach (\ref{eq_beta}) significantly reduces computational costs compared to the algorithm from  \cite{NgPer23} based on the entire sample size.
\end{document}